\documentclass[12pt, leqno]{amsart}

\usepackage{graphicx}
\usepackage{amssymb,amsmath,amsthm,amscd}
\usepackage{mathrsfs}
\usepackage{enumerate}
\usepackage{xspace}
\usepackage[usenames,dvipsnames]{color}
\usepackage[colorlinks=true, pdfstartview=FitV,
 linkcolor=blue,citecolor=blue,urlcolor=blue]{hyperref}
\usepackage[all]{xy}

\allowdisplaybreaks[4]

\newcounter{mycounter}

\newcounter{mycounters}

\newcounter{mycounterf}
\newlength{\my}
\newcommand{\nc}{\newcommand}
\numberwithin{equation}{section}

\newenvironment{nospace}[1][1pt]{\relax\setlength{\arraycolsep}{#1}
}{\relax}

\nc{\bns}{\begin{nospace}}
\nc{\ens}{\end{nospace}}

\newenvironment{myarray}[1]{\relax\setlength{\arraycolsep}{1pt}

\begin{array}{#1}}{\end{array}\relax}

\newcommand{\ba}{\begin{myarray}}
\newcommand{\ea}{\end{myarray}}

\nc{\bal}{\begin{myalign}}
\nc{\eal}{\end{myalign}}

\theoremstyle{plain}
\newtheorem{lemma}{Lemma}[section]
\newtheorem{prop}[lemma]{Proposition}
\newtheorem{theorem}[lemma]{Theorem}
\newcommand{\Prop}{\begin{prop}}
\newcommand{\enprop}{\end{prop}}
\newcommand{\Lemma}{\begin{lemma}}
\newcommand{\enlemma}{\end{lemma}}
\newcommand{\Th}{\begin{theorem}}
\newcommand{\enth}{\end{theorem}}
\newtheorem{corollary}[lemma]{Corollary}
\newcommand{\Cor}{\begin{corollary}}
\newcommand{\encor}{\end{corollary}}
\newtheorem{definition}[lemma]{Definition}

\newcommand{\Def}{\begin{definition}}
\newcommand{\edf}{\end{definition}}
\newtheorem{sublemma}[lemma]{Sublemma}
\newcommand{\Sublemma}{\begin{sublemma}}
\newcommand{\ensub}{\end{sublemma}}

\theoremstyle{definition}
\newtheorem{remark}[lemma]{Remark}
\newtheorem{example}[lemma]{Example}
\newtheorem*{Convention}{Convention}
\newcommand{\Conv}{\begin{Convention}}
\newcommand{\enconv}{\end{Convention}}
\nc{\Ex}{\begin{example}}
\nc{\enEx}{\end{example}}
\nc{\Rem}{\begin{remark}}
\nc{\enrem}{\end{remark}}

\newcommand{\C}{{\mathbb C}}

\newcommand{\Z}{{\mathbb Z}}
\newcommand{\B}{{\mathbf{B}}}

\newcommand{\R}{{\mathbb{R}}}

\newcommand{\codim}{{\operatorname{codim}}}
\newcommand{\Spec}{{\operatorname{Spec}}}

\newcommand{\seteq}{\mathbin{:=}}

\newcommand{\eim}[1]{{#1}_{\mspace{2mu}!}\mspace{1mu}}
\newcommand{\roim}[1]{\RR{#1}_*}
\newcommand{\reim}[1]{\RR{#1}_!}

\newcommand{\opb}[1]{#1^{-1}}

\newcommand{\epb}[1]{#1^{\,!}\,}

\newcommand{\g}{{\mathfrak{g}}}

\newcommand{\Hom}{\operatorname{Hom}}

\newcommand{\pt}{{\operatorname{pt}}}

\newenvironment{myequation}
{\relax\setlength{\arraycolsep}{1pt}\begin{eqnarray}}
{\end{eqnarray}}
\newenvironment{myequationn}
{\relax\setlength{\arraycolsep}{1pt}\begin{eqnarray*}}
{\end{eqnarray*}}

\nc{\eq}{\begin{myequation}}
\nc{\eneq}{\end{myequation}}
\nc{\eqn}{\begin{myequationn}}
\nc{\eneqn}{\end{myequationn}}

\newcommand{\hs}{\hspace*}
\nc{\ms}{\mspace}

\newcommand{\To}[1][{\hs{2ex}}]{\xrightarrow{\,#1\,}}

\newcommand{\on}{\operatorname}
\newcommand{\Ker}{\on{Ker}}

\newcommand{\bna}{\be[{\rm(a)}]}

\newcommand{\QED}{\end{proof}}
\newcommand{\Proof}{\begin{proof}}
\newcommand{\coh}{{\on{coh}}}
\newcommand{\soplus}{\mathop{\scalebox{0.6}{$\displaystyle\bigoplus$}}\limits}

\newcommand{\Supp}{\operatorname{Supp}}
\newcommand{\id}{\on{id}}

\newcommand{\Coker}{{\operatorname{Coker}}}

\newcommand{\bi}{\begin{enumerate}[{\rm(i)}]}

\newcommand{\set}[2]{\left\{#1 \mid #2 \right\}}

\newcommand{\Mod}{\operatorname{Mod}}

\newcommand{\eqsub}{\begin{subequations}\begin{eqnarray}}
\newcommand{\eneqsub}{\end{eqnarray}\end{subequations}}

\newcommand{\ol}{\overline}

\nc{\la}{\lambda}
\nc{\lam}{\lambda}
\nc{\U}[1][\g]{U_q(#1)}
\nc{\te}{\tilde{e}}
\nc{\tei}{\tilde{e}_i}
\nc{\tf}{\tilde{f}}
\nc{\tfi}{\tilde{f}_i}
\nc{\tU}{\widetilde U_q(\g)}
\nc{\tE}{\tilde{E}}
\nc{\tF}{\widetilde{\F}}
\nc{\tK}{\widetilde{K}}

\nc{\tk}{\tilde{k}}
\nc{\tkone}{\tk_{\ol{1}}}
\nc{\teone}{\tilde{e}_{\ol{1}}}
\nc{\tfone}{\tilde{f}_{\ol{1}}}

\nc{\teibar}{\tilde{e}_{\ol{i}}} \nc{\tfibar}{\tilde{f}_{\ol{i}}}
\nc{\tki}{{\tk}_{\ol {i}}}

\nc{\BZ}{{\mathbb{Z}}}
\nc{\al}{\alpha}
\nc{\qs}{{q}}
\nc{\lan}{\langle}
\nc{\ran}{\rangle}
\nc{\re}{{\mathrm{re}}}
\nc{\wt}{\operatorname{wt}}
\nc{\ch}{\operatorname{ch}}
\nc{\Uf}[1][\g]{U^-_q(#1)}
\nc{\Ue}{U^+_q(\g)}
\nc{\eps}{\varepsilon}
\nc{\vphi}{\varphi}
\nc{\sphi}{\varphi^*}
\nc{\seps}{\varepsilon^*}

\nc{\nn}{\nonumber}
\def\max{{\mathop{\mathrm{max}}}}
\nc{\vp}{\varpi}
\nc{\cls}{{\operatorname{cl}}}
\nc{\Wt}{{\operatorname{Wt}}}
\nc{\Us}{U'_q(\g)}
\nc{\La}{\Lambda}
\nc{\ro}{{\rm(}}
\nc{\rf}{{\rm)}\xspace}
\nc{\norm}{{\mathrm{norm}}}
\nc{\qbox}{\quad\mbox}
\nc{\braid}{{\mathfrak{B}}}
\nc{\Ad}{\operatorname{Ad}}
\nc{\Aut}{\operatorname{Aut}}

\nc{\Sn}{S^{{\mathrm{norm}}}}
\nc{\aff}{{\rm{aff}}}
\nc{\rk}{{\mathrm{rk}}}
\nc{\tP}{\widetilde{P}}
\nc{\tW}{\widetilde{W}}
\nc{\Dyn}{\mathrm{Dyn}}
\nc{\height}{{\operatorname{ht}}}
\nc{\bl}{\bigl(}
\nc{\br}{\bigr)}
\nc{\Hecke}{\mathrm{H}}
\nc{\HA}{\Hecke^{\mathrm{A}}}
\nc{\HB}{\Hecke^{\mathrm{B}}}
\newcommand{\scbul}{{\,\raise1pt\hbox{$\scriptscriptstyle\bullet$}\,}}
\nc{\vac}{{\phi}}
\nc{\Bt}{\B_\theta(\g)}
\nc{\be}{\begin{enumerate}}
\nc{\ee}{\end{enumerate}}
\nc{\low}{{\mathrm{low}}}
\nc{\upper}{{\mathrm{up}}}
\nc{\Zodd}{\Z_{\mathrm{odd}}}
\nc{\Ft}[1][n]{\mathbb{P}\mathrm{ol}_{#1}}
\nc{\Ftf}[1][n]{\widetilde{\mathbb{P}\mathrm{ol}}_{#1}}
\nc{\KA}{\on{K}^{\mathrm{A}}}
\nc{\KB}{\on{K}^{\mathrm{B}}}
\nc{\Res}{\on{Res}}
\nc{\Fc}[1][{n,m}]{\mathbf{F}_{#1}}
\nc{\tphi}{\tilde{\varphi}}
\nc{\CO}{\mathscr{O}}
\nc{\inte}{\mathrm{int}}
\nc{\Oint}{\mathcal{O}^{\ge0}_{\inte}}
\nc{\vs}{\vspace*}
\nc{\tL}{\widetilde{L}}
\nc{\tu}{\tilde{u}}
\nc{\noi}{\noindent}
\nc{\heigh}{\mathfrak{t}}
\nc{\lowest}{\mathfrak{l}}
\nc{\rootl}{\mathsf{Q}}
\nc{\cl}{\colon}
\nc{\uqpg}{U'_q(\mathfrak g)}
\nc{\Oh}{\widehat{\mathcal{O}}}

\nc{\KLR}{Khovanov-Lauda-Rouquier algebra}
\nc{\KLRs}{Khovanov-Lauda-Rouquier algebras}
\nc{\cor}{\mathbf{k}}
\nc{\cora}{{\cor(A)}}
\nc{\haut}{\mathrm{ht}}
\nc{\tens}{\mathop\otimes}
\nc{\gmod}{\mbox{-$\mathrm{gmod}$}}
\nc{\proj}{\mbox{-$\mathrm{proj}$}}
\nc{\gproj}{\mbox{-$\mathrm{gproj}$}}
\nc{\smod}{\mbox{-$\mathrm{mod}$}}

\nc{\h}{\mathfrak h}
\nc{\Rnorm}{R^{\rm{norm}}}
\nc{\K}{\C(q)}
\nc{\Vhat}{\widehat{V}}
\nc{\F}{\mathcal{F}}

\def\T{{\mathcal T}}

\nc{\fd}[1][A]{\on{\mathrm{flat.dim}_{#1}}}
\nc{\bP}{{\mathbb{P}}}
\nc{\bPh}{\widehat{\mathbb{P}}}
\nc{\bK}[1][{n}]{\widehat{\mathbb{K}}_{#1}}
\nc{\bV}[1][{n}]{\widehat{V}^{\otimes{#1}}}
\nc{\bVK}[1][{n}]{\widehat{V}^{\otimes{#1}}_{\widehat{\mathbb{K}}}}
\nc{\hV}{\widehat{V}}
\nc{\opp}{\mathrm{opp}}
\nc{\col}{\colon}
\nc{\bnum}{\be[{\rm(i)}]}
\nc{\oep}{\epsilon}
\nc{\qtext}{\quad\text}
\nc{\qtextq}[1]{\quad\text{#1}\quad}
\nc{\longtwoheadrightarrow}[1][]{\xymatrix{\ar@{->>}[r]^-{{#1}}&}}
\nc{\epiTo}[1][]{\longtwoheadrightarrow[{#1}]}
\nc{\epito}{\twoheadrightarrow}
\nc{\monoTo}[1][]{\xymatrix{\ar@{>->}[r]^-{{#1}}&}}
\nc{\sym}{\mathfrak{S}}
\nc{\inp}[1]{{({#1})_{\mathrm{n}}}}
\nc{\rtl}{\rootl}
\nc{\wtd}{\widetilde}
\nc{\etens}{\boxtimes}
\nc{\ds}[1]{\mathrm{d}(#1)}
\nc{\rmat}[1]{{r}_{\mspace{-2mu}\raisebox{-.5ex}{${\scriptstyle{#1}}$}}}
\nc{\shc}{\mathcal{C}}
\nc{\Fct}{{\on{Fct}}}
\nc{\tC}{\widetilde{\shc}}
\nc{\Zp}{\Z_{\ge0}}
\nc{\tPhi}{\widetilde{\Phi}}
\nc{\tT}{{\tilde{\T}}}
\nc{\Ob}{\on{Ob}}
\nc{\bwr}{\scalebox{1.3}{$\wr$}}
\nc{\Img}{\on{Im}}
\renewcommand{\Im}{\Img}
\nc{\Ab}{\mathcal{A}^{\mathrm{big}}}
\nc{\Sb}{\mathcal{S}^{\mathrm{big}}}
\nc{\As}{\mathcal{A}}
\nc{\ntens}{\widetilde{\otimes}}
\nc{\hR}{\widehat{R}}
\nc{\nconv}{\star}
\nc{\ts}{\tilde{s}}
\nc{\sho}{\mathscr{O}}
\nc{\bc}{\begin{cases}}
\nc{\ec}{\end{cases}}
\nc{\slnh}{{\widehat{\mathfrak{sl}}_N}}
\nc{\UA}{U_q'(\slnh)}
\nc{\KR}{R_K}
\nc{\cQ}{\mathcal{Q}}
\nc{\Irr}{\mathcal{I}rr}
\nc{\tQ}{\widetilde{\cQ}}
\nc{\bs}{\mathbf{s}}
\nc{\bL}{\mathbb{L}}
\nc{\tg}{\tilde{g}}

\newcommand{\Coim}{\on{Coim}}

\newcommand{\sect}{\on{\mathrm{\Gamma}}}

\newcommand{\rC}{\mathrm{C}}

\newcommand{\rD}{\mathrm{D}}

\newcommand{\Der}{\on{D}}
\newcommand{\Db}{\operatorname{D^{\mathrm{b}}}}

\newcommand{\sett}[2]{\set{#1}{\text{#2}}}

\nc{\Trg}{\mathcal{T}}
\nc{\scup}[1][]{\scalebox{0.8}{$\displaystyle\bigcup$}\mspace{2mu}
\raisebox{-.8ex}{$\scriptstyle #1$}\mspace{2mu}}
\nc{\scap}[1][]{\scalebox{0.8}{$\displaystyle\bigcap$}\mspace{2mu}
\raisebox{-.8ex}{$\scriptstyle #1$}\mspace{2mu}}
\nc{\ssqcup}[1][]%
{\scalebox{0.8}{$\displaystyle\bigsqcup$}\mspace{2mu}
\raisebox{-.8ex}{$\scriptstyle #1$}\mspace{2mu}}%
\nc{\ssqcupl}%
{\mathop{\scalebox{0.8}{$\displaystyle\bigsqcup$}}\limits}
\nc{\Tor}{\mathsf{T}}
\nc{\TF}{\mathsf{F}}
\nc{\cS}{\mathcal{S}}
\nc{\cons}{\mspace{.5mu}\mathrm{c}}
\nc{\Rcons}{{\R\text{-}\cons}}
\nc{\Ccons}{{\C\text{-}\cons}}
\nc{\pDc}{{}^{\mathrm{p}}\mspace{-3mu}\Der}
\nc{\pD}{{}^{1/2}\mspace{-3mu}\Der}
\nc{\ptau}{{}^{1/2}\mspace{0mu}\tau\mspace{2mu}}
\nc{\DX}[1][{\mathrm{b}}]{\Der^{{#1}}_\Rcons(\cor_X)}
\nc{\pDX}[1][{\mathrm{b}}]{\pD^{{#1}}_\Rcons(\cor_X)}
\nc{\tD}[1][X]{\mathsf{D}_{#1}}
\nc{\RR}{\mathrm{R}}
\nc{\rhom}[1][]{{\RR\mathscr{H}\mspace{-3mu}om}_{\raise1.5ex\hbox to.1em{}#1}}
\nc{\shom}[1][]{{\mathscr{H}\mspace{-3mu}om}_{\raise1.5ex\hbox to.1em{}#1}}
\nc{\rHom}[1][]{{\RR\mathrm{Hom}}_{\raise1.5ex\hbox to.1em{}#1}}
\nc{\rsect}{\mathrm{R}\Gamma}
\nc{\ori}{\mathrm{or}}
\nc{\DA}[1][\mathrm{b}]{\Der^{#1}_{\coh}(A)}
\nc{\pDA}[1][\mathrm{b}]{\pD^{#1}_{\coh}(A)}
\nc{\DO}[1][\mathrm{b}]{\Der^{#1}_{\coh}(\sho_X)}
\nc{\pDO}[1][\mathrm{b}]{\pD^{#1}_{\coh}(\sho_X)}
\nc{\shf}{\mathscr{F}}
\nc{\shh}{\mathscr{H}}
\nc{\shg}{\mathscr{G}}
\nc{\shd}{\mathcal{D}}
\nc{\cF}{\mathcal{F}}
\nc{\DXA}[1][{\mathrm{b}}]{\Der^{#1}_{\Rcons}(A_X)}
\nc{\DXAC}[1][{\mathrm{b}}]{\Der^{#1}_{\Ccons}(A_X)}
\nc{\DYA}[1][{\mathrm{b}}]{\Der^{#1}_{\Rcons}(A_Y)}
\nc{\pDXA}[1][{\mathrm{b}}]{\pD^{#1}_{\Rcons}(A_X)}
\nc{\pDXAC}[1][{\mathrm{b}}]{\pD^{#1}_{\Ccons}(A_X)}
\nc{\pDYA}[1][{\mathrm{b}}]{\pD^{#1}_{\Rcons}(A_Y)}
\nc{\pDXT}[1][{\mathrm{b}}]{\pD^{#1}_{\Ccons}(A_{T^*X})}
\nc{\pdim}{\on{mod\text{-}dim}}
\nc{\dX}{{\overset{\circ}{X}}}
\nc{\di}{{\widetilde{i}}}
\nc{\good}{good\xspace}
\newcommand{\RHom}[1][]{\RR\mathrm{Hom}_{\raise1.5ex\hbox to.1em{}#1}}
\nc{\dt}{distinguished triangle\xspace}
\nc{\dts}{distinguished triangles\xspace}
\nc{\tone}{\To[\;+1\;]}
\nc{\ks}{{}^{\ms{-2mu}1/2}_{\mathrm{KS}}\mspace{-1mu}\Der}
\nc{\KS}[1][{\mathrm{b}}]{\ks^{#1}_{\Rcons}(A_X)}
\nc{\KSC}[1][{\mathrm{b}}]{\ks^{#1}_{\Ccons}(A_X)}
\nc{\dc}{\mathop{\dim_{\mspace{2mu}\C}}}
\nc{\codc}{\mathop{\mathrm{codim}_{\mspace{2mu}\C}}}

\nc{\KST}[1][{\mathrm{b}}]{\ks^{#1}_{\Ccons}(A_{T^*X})}
\nc{\KSXY}[1][{\mathrm{b}}]{\ks^{#1}_{\Ccons}(A_{T^*_YX})}

\newcommand{\letens}{\overset{\mathrm{L}}{\etens}}

\nc{\muhom}{\mu hom}
\nc{\pDT}[1][{\mathrm{b}}]{\pD^{#1}_{\Ccons}(A_{T^*X})}
\newcommand{\ltens}[1][]{\mathbin{\overset{\mathrm{L}}\tens}_{#1}}
\nc{\DXCk}{\Der^{\mathrm{b}}_{\Ccons}(\cor_X)}
\nc{\pDXCk}[1][{\mathrm{b}}]{\pDc^{#1}_{\Ccons}(\cor_X)}
\nc{\DXk}{\Der^{\mathrm{b}}_{\Rcons}(\cor_X)}
\nc{\pDXk}[1][{\mathrm{b}}]{\pD^{#1}_{\Rcons}(\cor_X)}
\nc{\nos}[1][0pt]{\setlength{\arraycolsep}{#1}}
\nc{\Ss}{\on{SS}}
\nc{\tstr}[2]{\bl#1, #2\br_{c\in\R}}
\renewcommand{\le}{\leqslant}
\renewcommand{\ge}{\geqslant}
\renewcommand{\Re}{\mathrm{Re}}

\newlength{\mylength}
\begin{document}

\title{Self-dual T-structure}

\author{Masaki Kashiwara}
\thanks{The research was supported by Grant-in-Aid for Scientific Research (B) 
15H03608, Japan Society for the Promotion of Science.}
\address{Research Institute for Mathematical Sciences \\
          Kyoto University \\ Kyoto 606-8502, Japan}

\email{masaki@kurims.kyoto-u.ac.jp}


\begin{abstract}
We give a self-dual t-structure on the derived category
of $\mathbb{R}$-constructible sheaves 
over any Noetherian regular ring
by generalizing the notion of t-structure.
\end{abstract}

\date{}

\maketitle

\section*{Introduction}
Let $X$ be a complex manifold and let $\DXCk$ be the derived category of
sheaves of $\cor$-vector spaces on $X$ with $\C$-constructible cohomologies.
Here $\cor$ is a given base field.
Then the t-structure $\bl\pDXCk[\le0],\pDXCk[\ge0]\br$
on $\DXCk$ with  middle perversity is 
self-dual with respect to the Verdier dual functor 
$\tD=\rhom(\scbul,\omega_X)$. Namely, the Verdier dual functor
exchanges $\pDXCk[\le0]$ and $\pDXCk[\ge0]$.
However, on a real analytic manifold $X$ (of positive dimension), 
any perversity does not give a self-dual
t-structure on
the derived category $\DXk$ of $\R$-constructible sheaves on $X$.
In this paper, we construct such a self-dual t-structure after generalizing the notion of t-structure.
This generalized notion  already appeared in the paper of Bridgeland on
stability conditions (\cite{B07}). (See also \cite{K08}.)
This construction can be also applied to the derived category $\DA$ of
finitely generated modules over a Noetherian regular ring $A$.
We construct a (generalized) t-structure on $\DA$ which is
self-dual with respect to the duality functor
$\RHom[A](\scbul,A)$.

Let us explain our results more precisely in the example of
$\DXk$. Let $X$ be a real analytic manifold.
Recall that a sheaf $F$
of $\cor$-vector space is called $\R$-constructible if 
$X$ is a locally finite union of locally closed subanalytic subsets $\{X_\al\}_\al$ 
such that all the restrictions $F\vert_{X_\al}$ are locally constant with finite-dimensional fibers.
Let $\DXk$ be the bounded derived category of $\R$-constructible sheaves.
Let $\tD=\rhom(\scbul,\omega_X)$ be the Verdier dual functor.
For $c\in\R$, we define
\eq
&&\hs{4ex}\ba{l}\pDXk[\le c]\seteq\set{K\in\DXk}{
\text{$\dim\Supp(H^iK)\le 2(c-i)$ for any $i\in\Z$}},\\[1ex]
\pDXk[\ge c]\seteq\set{K\in\DXk}{\tD K\in\pDXk[\le -c]}.
\ea
\eneq
Then, the pair $\bl(\pDXk[\le c])_{c\in\R},(\pDXk[\ge c])_{c\in\R}\br$ satisfies the
axiom of (generalized) t-structure (Definition~\ref{axiom}).
In particular, $\bl\pDXk[\le c],\pDXk[> c-1]\br$ 
is a t-structure in the ordinary sense for any $c\in\R$. Here 
$\pDXk[>c]\seteq\scup[{b>c}]\pDXk[\ge b]$.
Therefore,
for any $K\in\DXk$ and $c\in\R$, there exists a \dt
$K'\to K\to K''\tone$
in $\DXk$ such that $K'\in\pDXk[\le c]$
and  $K''\in\pDXk[>c]$.

Note that $\pDXk[\le c]=\pDXk[\le s]$ for $s\in\frac{1}{2}\Z$ such that
$s\le c<s+1/2$, and
$\pDXk[\ge c]=\pDXk[\ge s]$ for $s\in\frac{1}{2}\Z$ such that
$s-1/2<c\le s$. 

\medskip
This paper is organized as follows.
In Section \ref{sec:abel}, we generalize the notion of a t-structure.
In Section \ref{sec:2}, we recall the result of \cite{K08}
on a t-structure on the derived category of a quasi-abelian category.
In Section \ref{sec:3}, we give the t-structure associated with a torsion pair on an abelian category.

In Section \ref{sec:4}, we give a self-dual t-structure on the derived category of
coherent sheaves on a Noetherian regular scheme.

In Section \ref{sec:5}, we define two t-structures on the derived category
of the abelian category of $\R$-constructible sheaves of  $A$-modules
on a subanalytic space $X$. Here $A$ is a Noetherian regular ring.
One is purely topological and the other is self-dual
with respect to the Verdier duality functor.

In Section \ref{sec:6}, we study the self-dual t-structure on
the derived category
of the abelian category of sheaves of $A$-modules on a complex manifold $X$
with $\C$-constructible cohomologies. The main result is its 
microlocal characterization (Theorem \ref{th:micro}).

\begin{Convention}
In this paper, all subanalytic spaces and complex analytic spaces 
are assumed to be Hausdorff, locally compact, countable at infinity and 
with finite dimension.
\end{Convention}

\section{(Generalized) T-structure}
\label{sec:abel}

Since the following lemma is elementary, we omit its proof.
\Lemma Let $X$ be a set.
\bnum
\item Let $(X^{\le c})_{c\in\R}$ be a family of subsets of $X$
such that $X^{\le c}=\scap[{b>c}]X^{\le b}$ for any $c\in\R$.
Set $X^{<c}\seteq \scup[{b<c}]X^{\le b}$.
Then we have
\bna
\item $X^{<c}=\scup[{b<c}]X^{<b}$,
\item $X^{\le c}=\scap[{b>c}]X^{<b}$.
\ee
\item Conversely, let $(X^{<c})_{c\in\R}$ be a family of subsets of $X$
such that $X^{<c}=\scup[{b<c}]X^{<b}$ for any $c\in\R$.
Set $X^{\le c}\seteq \scap[{b>c}]X^{<b}$.
Then we have
\bna
\item $X^{\le c}=\scap[{b>c}]X^{\le b}$,
\item $X^{<c}=\scup[{b<c}]X^{\le b}$.
\ee
\item 
Let $(X^{\le c})_{c\in\R}$ and $(X^{<c})_{c\in\R}$ be as above.
Let $a,b\in \R$ such that $a<b$.
If $X^{<c}=X^{\le c}$ for any $c$ such that $a<c\le b$, then
$X^{\le a}=X^{\le b}$.
\ee
\enlemma

Let us recall the notion of t-structure (see \cite{BBDG81}).
Let $\Trg$ be a triangulated category.
Let $\Trg^{\le0}$ and $\Trg^{\ge0}$ be strictly full subcategories of $\Trg$.
Here, a subcategory $\shc'$ of a category $\shc$ 
is called {\em strictly full} if
it is full, i.e. $\Hom_{\shc'}(X,Y)=\Hom_\shc(X,Y)$ for any $X,Y\in\shc'$,
and any object of $\shc$ 
isomorphic to some object of $\shc'$ is an object of $\shc'$.

For $n\in\Z$, we set
$\Trg^{\le n}=\Trg^{\le0}[-n]$ and $\Trg^{\ge n}=\Trg^{\ge0}[-n]$.
Let us recall that $(\Trg^{\le0},\Trg^{\ge0})$ is a t-structure on $\Trg$
 if it satisfies:
\eq&&\left\{\parbox{70ex}{
\bna
\item $\Trg^{\le0}\subset\Trg^{\le1}$ and 
$\Trg^{\ge1}\subset\Trg^{\ge0}$,
\item $\Hom_{\Trg}(X,Y)=0$ for $X\in\Trg^{\le0}$ and $Y\in\Trg^{\ge1}$,
\item for any $X\in\Trg$, there exists a \dt
$X_0\to X\to X_1\tone$ in $\Trg$ such that $X_0\in \Trg^{\le0}$ and
$X_1\in\Trg^{\ge1}$.
\ee
}\right.\eneq

We shall generalize this notion. 
\Def\label{axiom}
Let
$(\Trg^{\le c})_{\;c\in\R}$ and $(\Trg^{\ge c})_{\;c\in\R}$ 
be families of strictly full subcategories of a triangulated category $\Trg$,
and set
$\Trg^{< c}=\scup[{b<c}]\Trg^{\le b}$ and $\Trg^{> c}=\scup[{b>c}]\Trg^{\ge b}$.
We say that $\bl (\Trg^{\le c})_{c\in\R},(\Trg^{\ge c})_{c\in\R}\br$ 
is a 
{\em \ro generalized\rf\ t-structure} \ro cf.\ \cite{B07}\rf\ if it satisfies
the following conditions.
\eq&&\left\{\parbox{70ex}{
\bna
\item $\Trg^{\le c}=\scap[{b>c}]\Trg^{\le b}$ and 
$\Trg^{\ge c}=\scap[{b<c}]\Trg^{\ge b}$ for any $c\in\R$,\label{ax:1}
\item $\Trg^{\le c+1}=\Trg^{\le c}[-1]$ and $\Trg^{\ge c+1}=\Trg^{\ge c}[-1]$ for any $c\in\R$,\label{ax:2}
\item
we have $\Hom_{\Trg}(X,Y)=0$ for any $c\in\R$,
$X\in\Trg^{<c}$ and  $Y\in\Trg^{>c}$,\label{ax:3}
\item  
for any $X\in\Trg$ and $c\in\R$, there exist \dts
$X_0\to X\to X_1\tone$ and $X'_0\to X\to X'_1\tone$ 
in $\Trg$ such that $X_0\in \Trg^{\le  c}$, $X_1\in\Trg^{ > c}$
and $X'_0\in \Trg^{<c}$, $X'_1\in\Trg^{\ge c}$.\label{ax:4}
\setcounter{mycounter}{\value{enumi}}
\ee
}\right.\eneq
\edf

Note that under conditions (a)--(c), the \dts in (d) are unique 
up to a unique isomorphism.

\smallskip
If $\bl (\Trg^{\le c})_{c\in\R},(\Trg^{\ge c})_{c\in\R}\br$ is a generalized t-structure,
then the pairs $(\Trg^{\le c},\Trg^{>c-1})$ and  $(\Trg^{< c},\Trg^{\ge c-1})$ are t-structures
in the original sense for any $c\in\R$.
Hence, $\Trg^{\le c}\cap\Trg^{>c-1}$ and $\Trg^{<c}\cap\Trg^{\ge c-1}$
are abelian categories.

Assume that $\bl (\Trg^{\le c})_{c\in\R},(\Trg^{\ge c})_{c\in\R}\br$ is a generalized t-structure.
Then the inclusion functors
$\Trg^{\le c}\To\Trg$ and $\Trg^{< c}\To\Trg$ have right adjoints
$$\text{
$\tau^{\le c}\cl \Trg\To \Trg^{\le c}$ and
$\tau^{<c}\cl \Trg\To \Trg^{<c}$, respectively.}$$
Similarly, the inclusion functors
$\Trg^{\ge c}\To\Trg$ and $\Trg^{> c}\To\Trg$ have left adjoints
$$\text{
$\tau^{\ge c}\cl \Trg\To \Trg^{\ge c}$ and
$\tau^{>c}\cl \Trg\To \Trg^{>c}$, respectively.}$$
We have  \dts functorially in $X\in\Trg$:
\eqn
&&\tau^{\le c}X\To X\To\tau^{>c}X\tone\quad\text{and}\\
&&\tau^{<c}X\To X\To\tau^{\ge c}X\tone.
\eneqn
These four functors are called the truncation functors of the generalized t-structure
$\bl (\Trg^{\le c})_{c\in\R},(\Trg^{\ge c})_{c\in\R}\br$.

For any $a,b\in\R$, we have isomorphisms of functors:
\eqn
&&\tau^{\le a}\circ\tau^{\le b}\simeq\tau^{\le \min(a,b)}
\,,\quad
\tau^{\ge a}\circ\tau^{\ge b}\simeq\tau^{\ge \max(a,b)}
\ \text{, and}\\
&& \tau^{\le a}\circ\tau^{\ge b}\simeq\tau^{\ge b}\circ\tau^{\le a}.
\eneqn
In the last formula, we can replace $\tau^{\ge a}$ with $\tau^{>a}$
or $\tau^{\le b}$ with $\tau^{<b}$.
For any $c\in\R$, we have
\eq
&&\ba{rl}
\Trg^{\le c}&=\set{X\in\Trg}%
{\text{$\Hom_\Trg(X,Y)\simeq 0$ for any $Y\in\Trg^{>c}$}},\\
\Trg^{<c}&=\set{X\in\Trg}%
{\text{$\Hom_\Trg(X,Y)\simeq 0$ for any $Y\in\Trg^{\ge c}$}},\\
\Trg^{\ge c}&=\set{Y\in\Trg}%
{\text{$\Hom_\Trg(X,Y)\simeq 0$ for any $X\in\Trg^{<c}$}},\\
\Trg^{>c}&=\set{Y\in\Trg}%
{\text{$\Hom_\Trg(X,Y)\simeq 0$ for any $X\in\Trg^{\le c}$}}.
\ea\label{eq:orth}
\eneq

We set $\Trg^c\seteq\Trg^{\le c}\cap\Trg^{\ge c}$. 
Then $\Trg^c$ is a quasi-abelian category
(see \cite{B07} and \cite{Sch04}).
More generally, for $a\le b$, we set 
$$\Trg^{[a,b]}\seteq\Trg^{\le b}\cap\Trg^{\ge a}.$$
Then $\Trg^{[a,b]}$ is a quasi-abelian category if $a\le b<a+1$.

A t-structure $(\Trg^{\le0},\Trg^{\ge0})$ is regarded 
as a generalized t-structure by
\eq&&\ba{rll}
\Trg^{\le c}&=\Trg^{\le 0}[-n]\quad
&\text{for $n\in\Z$ such that $n\le c<n+1$,}\\
\Trg^{\ge c}&=\Trg^{\ge 0}[-n]\quad
&\text{for $n\in\Z$ such that $n-1< c\le n$.}
\ea
\eneq
Hence, a t-structure is nothing but a generalized t-structure such that
$\Trg^{\le0}=\Trg^{<1}$ and $\Trg^{\ge1}=\Trg^{>0}$,
or equivalently $\Trg^c=0$ for any $c\notin \Z$.

\smallskip
{\em In the sequel, we call
a generalized t-structure simply a t-structure.}

\Rem
In the examples we give in this paper, the t-structures 
also satisfy the following condition:

\bna
\setcounter{enumi}{\value{mycounter}}
\item
for any $c\in \R$ we can find $a$ and $b$ such that  $a<c<b$ and
\be[(1)]
\item $\Trg^{<c}=\Trg^{\le a}$, $\Trg^{\le c}=\Trg^{<b}$,
\item $\Trg^{>c}=\Trg^{\ge b}$, $\Trg^{\ge c}=\Trg^{>a}$.
\ee
\ee
More precisely, in the examples in this paper,
we can take $a=\max\set{s\in\frac{1}{2}\Z}{s<c}$ and
$b=\min\set{s\in\frac{1}{2}\Z}{s>c}$. Hence $\Trg^c=0$ if $c\not\in\frac{1}{2}\Z$.
\enrem

\section{T-structure on the derived category of a quasi-abelian category}\label{sec:2}
For more details, see \cite[{\S\;2}]{K08}. 

Let $\shc$ be a quasi-abelian category
(see \cite{Sch04}).
Recall that, for a morphism $f\cl X\to Y$ in $\shc$,
$\Im(f)\seteq\Ker\bl Y\to\Coker(f)\br$ and $\Coim(f)\seteq\Coker\bl \Ker(f)\to X\br$. Hence, we have a diagram
$$\xymatrix
{\Ker(f)\ar[r]&X\ar[r]\ar@/^1.5pc/[rrr]|f&\Coim(f)\ar[r]&\Im(f)\ar[r]&Y\ar[r]&\Coker(f)}
.$$
Let $\rC(\shc)$ be the category of complexes in $\shc$, and
$\rD(\shc)$ the derived category of $\shc$ (see \cite{Sch04}).
Let us define various truncation functors 
\index{truncation}%
for $X\in\rC(\shc)$:
\eqn\hs{3ex}
\ba{rcccccccccccccccc}
\tau^{\le n}X&\hs{1ex}:\hs{2ex}& \cdots&\to&X^{n-1}
&\to&\Ker d_X^n&\to&0&\to&0&\to&\cdots\\[1ex]
\tau^{\le n+1/2}X&:& \cdots&\to& X^{n-1}
&\to& X^n&\to&\Im d_X^n&\to&0&\to&\cdots\\[1ex]
\tau^{\ge n}X&:& \cdots&\to&0
&\to&\Coker d_X^{n-1}&\to& X^{n+1}&\to &X^{n+2}&\to&\cdots\\[1ex]
\tau^{\ge n+1/2}X&:& \cdots&\to&0&\to&\Coim d_X^{n}&\to& X^{n+1}
&\to& X^{n+2}&\to&\cdots
\ea
\eneqn
for $n\in\Z$.
Then we have morphisms functorial in $X$:
$$\tau^{\le s}X\To \tau^{\le t}X\To X\To\tau^{\ge s}X\To
\tau^{\ge t}X$$
for $s,t\in\frac{1}{2}\Z$ such that $s\le t$.
We can easily check that these functors
$\tau^{\le s},\tau^{\ge s}\cl\rC(\shc)\to\rC(\shc)$
send the morphisms homotopic to zero to
morphisms homotopic to zero and the quasi-isomorphisms to
quasi-isomorphisms.
Hence, they induce the functors
$$\tau^{\le s},\tau^{\ge s}\cl\rD(\shc)\to\rD(\shc)$$
and morphisms $\tau^{\le s}\to\id\to\tau^{\ge s}$.

For $s\in\frac{1}{2}\Z$, set
\eqn
\rD^{\le s}(\shc)&=&\sett{X\in\rD(\shc)}%
{$\tau^{\le s}X\to X$ is an isomorphism}\\
\rD^{\ge s}(\shc)&=&\sett{X\in\rD(\shc)}%
{$X\to\tau^{\ge s}X$ is an isomorphism}.
\eneqn

Then $\{\rD^{\le s}(\shc)\}_{s\in\frac{1}{2}\Z}$
is an increasing sequence of strictly full subcategories of
$\rD(\shc)$, and $\{\rD^{\ge s}(\shc)\}_{s\in\frac{1}{2}\Z}$
is a decreasing sequence of strictly full subcategories of
$\rD(\shc)$.

The functor $\tau^{\le s}\cl\rD(\shc)\to\rD^{\le s}(\shc)$
is a right adjoint functor of the inclusion functor
$\rD^{\le s}(\shc)\hookrightarrow \rD(\shc)$,
and $\tau^{\ge s}\cl\rD(\shc)\to\rD^{\ge s}(\shc)$
is a left adjoint functor of $\rD^{\ge s}(\shc)\hookrightarrow\rD(\shc)$.

For $c\in\R$, we set
\eq
&&\ba{l}
\rD^{\le c}(\shc)=\rD^{\le s}(\shc)\quad\text{where
$s\in\frac{1}{2}\Z$ satisfies $s\le c<s+1/2$,}\\[1ex]
\rD^{\ge c}(\shc)=\rD^{\ge s}(\shc)\quad\text{where
$s\in\frac{1}{2}\Z$ satisfies $s-1/2< c\le s$.}
\ea
\label{eq:interp}
\eneq
Then
\Prop[{\cite{Sch04}, see also \cite{K08}}]
$\bl (\rD^{\le c}(\shc))_{c\in\R},(\rD^{\ge c}(\shc))_{c\in\C}\br$
is a t-structure.
\enprop

We call it the {\em standard t-structure} on $\rD(\shc)$.
The triangulated category $\rD(\shc)$ is equivalent to the
derived category of the
abelian category $\rD^{\le c}(\shc)\cap\rD^{> c-1}(\shc)$ for every $c\in\R$.
The full subcategory $\rD^0(\shc)\seteq\rD^{\le 0}(\shc)\cap\rD^{\ge 0}(\shc)$ 
is equivalent to $\shc$.

If $\shc$ is an abelian category, then the standard t-structure is:
\eqn
&&\rD^{\le c}(\shc)=\set{X\in\rD(\shc)}{\text{$H^i(X)=0$ for any $i>c$}},\\
&&\rD^{\ge c}(\shc)=\set{X\in\rD(\shc)}{\text{$H^i(X)=0$ for any $i<c$}}.
\eneqn

\section{T-structure associated with a torsion pair}\label{sec:3}

Let $\shc$ be an abelian category.
A torsion pair
is a pair $(\Tor,\TF)$ of strictly full subcategories of $\shc$
such that
\eq
&&\left\{\parbox{70ex}{
\bna
\item $\Hom_\shc(X,Y)=0$ for any $X\in\Tor$ and $Y\in\TF$,\vs{1ex}
\item for any $X\in\shc$, there exists an exact sequence
$0\to X'\to X\to X''\to 0$
with  $X'\in\Tor$ and $X''\in\TF$.
\ee
}\right.
\eneq
Let $(\Tor,\TF)$ be a torsion pair. Then we have
\eqn
\Tor&\simeq&\set{X\in\shc}{\text{$\Hom_\shc(X,Y)=0$ for any $Y\in\TF$}},\\
\TF&\simeq&\set{Y\in\shc}{\text{$\Hom_\shc(X,Y)=0$ for any $X\in\Tor$}}.
\eneqn
Moreover, $\Tor$ is stable under taking quotients and extensions,
while $\TF$ is stable under taking subobjects and extensions.

For any integer $n$, we define
\eq&&\hs{5ex}\ba{rl}
\pDc^{\le n}(\shc)&\seteq\set{X\in\Der(\shc)}{\text{$H^i(X)\simeq0$
for any $i>n$}},\\
\pDc^{\le n-1/2}(\shc)&\seteq\set{X\in\Der(\shc)}{\text{
$H^i(X)\simeq0$ for any $i>n$ and $H^n(X)\in\Tor$}},\\
\pDc^{\ge n-1/2}(\shc)&\seteq\set{X\in\Der(\shc)}{\text{$H^i(X)\simeq0$
for any $i<n$}},\\
\pDc^{\ge n}(\shc)&\seteq\set{X\in\Der(\shc)^{\ge n-1/2}}{\text{
$H^i(X)\simeq0$ for any $i<n$ and $H^n(X)\in\TF$}}.
\ea
\eneq

For any $c\in\R$, we define $\pDc^{\le c}(\shc)$
and $\pDc^{\ge c}(\shc)$ by \eqref{eq:interp}.\\
Since the following proposition can be easily proved, we omit the proof.
\Prop 
$\bl(\pDc^{\le c}(\shc))_{c\in\R},(\pDc^{\ge c}(\shc))_{c\in\R}\br$
is a t-structure.
\enprop
We have
\eqn
&&\Tor\simeq\pDc^{-1/2}(\shc),\quad\TF\simeq\pDc^0(\shc),
\qtextq{and} \shc\simeq\pDc^{[-1/2,0]}(\shc).
\eneqn
Moreover, $\Der(\shc)$ is equivalent to the derived category of
the abelian category $\pDc^{[0,1/2]}(\shc)$.

Note that we have
$$\Der^{\le c}(\shc)\subset \pDc^{\le c}(\shc)\subset\Der^{\le c+1/2}(\shc)\quad\text{and}\quad
\Der^{\ge c+1/2}(\shc)\subset\pDc^{\ge c}(\shc)\subset\Der^{\ge c}(\shc).$$

\section{Self-dual t-structure on 
the derived category of coherent sheaves}\label{sec:4}
Let $X$ be a Noetherian regular scheme.
Let $\tD[X]$ be the duality functor
$\tD[X]\seteq\rhom[{\sho_X}](\scbul,\sho_X)$.
Let $\DO$ be the bounded derived category of $\sho_X$-modules
with coherent cohomologies.
We denote by
$\bl(\DO[\le c])_{c\in\R},(\DO[\ge c])_{c\in\R}\br$ the standard t-structure on
$\DO$.

Recall that, for any coherent $\sho_X$-module $\shf$, its codimension is defined by
$$\codim(\shf)\seteq\codim\bl\Supp(\shf)\br=\inf_{x\in\Supp(\shf)}\dim\sho_{X,\,x}.$$ Here we understand $\codim(0)=+\infty$.

We set
\eqn
\pDO[\le c]&\seteq&\set{\shf\in\DO}%
{\text{$\codim(H^i(\shf))\ge 2(i-c)$ for any $i\in\Z$}},\\[1ex]
\pDO[\ge c]&\seteq&\set{\shf\in\DO}{\tD[X]\shf\in\DO[\le -c]}\\
&=&\set{\shf\in\DO}%
{\text{$\codim(H^i(\tD[X]\shf))\ge 2(i+c)$ for any $i\in\Z$}}.
\eneqn
Then they satisfy Definition~\ref{axiom} \eqref{ax:1}.
Remark that we have 
\eqn
&&\pDO[\le c]=\set{\shf\in\DO}%
{\text{$\shf_x\in\Der^{\le c+\frac{1}{2}\dim \sho_{X,x}}(\sho_{X,x})$
for any $x\in X$}}.
\eneqn
We have also
\eqn
\pDO[< c]&\seteq&\scup[{b<c}]\pDO[\le b]\\
&=&\set{\shf\in\DO}%
{\text{$\codim(H^i(\shf))>2(i-c)$ for any $i\in\Z$}},\\[1ex]
\pDO[>c]&\seteq&\scup[{b>c}]\pDO[\ge b]\\
&=&\set{\shf\in\DO}{\text{$\codim(H^i(\tD\shf))>2(i+c)$ for any $i\in\Z$}}.
\eneqn

\Lemma\label{lem:Ot}
Let $\shf\in\DO$.
Then $\shf\in \pDO[\ge c]$ if and only if we have
$H^i\rsect_Z\shf=0$ for any closed 
subset $Z$ and $i<c+\codim Z/2$.
\enlemma 

\Proof
We shall use the results in \cite{K04}.
Let us define the systems of support
\eqn
&&\Phi^n=\set{Z}{\codim Z\ge 2(n+c)},\\
&&\Psi^n=\set{Z}{n<c+1+\codim Z/2}.\\
\eneqn
Then it is enough to show that
\eq
&&(\Phi\circ \Psi)^n\seteq\scup[{i+j=n}]\bl\Phi^i\cap\Psi^j\br
=\set{Z}{\codim Z\ge n}.\label{eq:PhiPsi}
\eneq
Indeed, one has
$$\pDO[\le -c]={}^\Phi\mspace{-2mu}\DO[{\le 0}]
\seteq\set{\shf\in \DO}%
{\text{$\Supp(H^k(\shf))\in \Phi^k$ for any $k\in\Z$}}$$
 and hence
\cite[Theorem 5.9]{K04} along with \eqref{eq:PhiPsi}
implies that $\pDO[\ge c]$ coincides with
\eqn
{}^\Psi\mspace{-2mu}\DO[{\ge 0}]&\seteq&\set{F}%
{\text{$H^i(\rsect_ZF)=0$ for any $Z\in\Psi^{i+1}$}}\\
&=&\set{F}%
{\text{$H^i(\rsect_ZF)=0$ for any $i<c+\codim Z/2$}}.
\eneqn

Let us show \eqref{eq:PhiPsi}
Assume that
$Z\in \Phi^i\cap\Psi^j$ with $i+j=n$.
Then we have
$$2\,\codim Z\ge 2(i+c)+\bl 2(j-c-1)+1\br=2n-1$$
 and hence
$\codim Z\ge n$.

Conversely assume that $\codim Z\ge n$.
Then take an integer $i$ such that $i\le \codim Z/2-c<i+1$.
Then we have $i>\codim Z/2-c-1$ and
$$j\seteq n-i<\codim Z-(\codim Z/2-c-1)=c+1+\codim Z/2.$$
Hence $Z\in\Phi^i\cap\Psi^j\subset(\Phi\circ\Psi)^n$.
\QED

\Prop
$\bl(\pDO[\le c])_{c\in\R},(\pDO[\ge c])_{c\in\R}\br$ is a t-structure on
$\DO$.
\enprop

\Proof It follows from \cite{K04}. Indeed, 
the pair $\bl\pDO[<c+1],\pDO[\ge c]\br$
coincides with $\bl {}^\Psi\mspace{-2mu}\DO^{\le 0}, {}^\Psi\mspace{-2mu}\DO^{\ge 0}\br$
by the proof of the preceding proposition.
\QED

\Cor\label{cor:hom} For $\shf\in\pDO[\le c]$ and $\shg\in\pDO[\ge c']$, we have
$$\rhom[\sho_X](\shf,\shg)\in\DO[\ge c'-c].$$
Conversely we have
\eqn
\pDO[\ge c']&=&\{\shg\in\DO\mid
\parbox[t]{50ex}{$\rhom[\sho_X](\shf,\shg)\in\DO[\ge c'-c]$ for any 
$c\in\R$ and $\shf\in\pDO[\le c]$\}\quad \hfill for any $c'\in\R$,}\\[1ex]
\pDO[\ge c]&=&\{\shf\in\DO\mid
\parbox[t]{50ex}{$\rhom[\sho_X](\shf,\shg)\in\DO[\ge c'-c]$ for any 
$c'\in\R$ and $\shg\in\pDO[\ge c']$\}\quad\hfill for any $c\in\R$.}
\eneqn
\encor

\Prop\label{prop:inner}
For $\shf$, $\shg\in\DO$, we have
\bnum
\item
if $\shf\in\pDO[\le c]$ and $\shg\in\DO[\le c']$, 
then $\shf\ltens[{\sho_X}]\shg\in\pDO[\le c+c']$,

\vs{.5ex}
\item if $\shf\in\DO[\le c]$ and $\shg\in\pDO[\ge c']$, 
then $\rhom[\sho_X](\shf,\shg)\in\pDO[\ge c'-c]$,

\vs{.5ex}
\item
if $\shf\in\pDO[\ge c]$ and $\shg\in\DO[\le c']$, 
then $\rhom[\sho_X](\shf,\shg)\in\pDO[\le c'-c]$,

\item if $\shf\in\pDO[\ge c]$ and $\shg\in\pDO[\ge c']$, 
then $\shf\ltens[\sho_X]\shg\in\DO[\ge c+c']$.\label{item:3}
\ee
\enprop
\Proof
(i)\ For any $\shh\in \pDO[\ge c'']$,
we have $\rhom[\sho_X](\shf,\shh)\in\DO[\ge c''-c]$ by Corollary~\ref{cor:hom}.
Hence,
$\rhom[\sho_X](\shf\ltens[{\sho_X}]\shg,\shh)\simeq
\rhom[\sho_X]\bl \shg,\rhom[\sho_X](\shf,\shh)\br$
belongs to $\DO[\ge c''-c-c']$.
Since it holds for an arbitrary $\shh\in \pDO[\ge c'']$,
we conclude $\shf\ltens[{\sho_X}]\shg\in\pDO[\le c+c']$
by \eqref{eq:orth}.

\medskip\noi
(ii)\ 
Since $\shf\ltens[\sho_X]\tD\shg\in\pDO[\le c-c']$ by (i),
and hence
$\rhom[\sho_X](\shf,\shg)\simeq\tD\bl \shf\ltens\tD\shg\br$
belongs to $\pDO[\ge c'-c]$.

\medskip\noi
(iii)\ 
Since $\rhom[\sho_X](\shf,\shg)
\simeq(\tD\shf)\ltens[\sho_X]\shg$,
(iii) follows from (i).

\medskip\noi
(iv) follows from Corollary~\ref{cor:hom} and
$\shf\ltens[\sho_X]\shg\simeq\rhom[\sho_X](\tD \shf,\shg)$.
\QED
Let $A$ be a Noetherian regular ring and $X=\Spec(A)$.
We write $\DA$, $\pDA[\le c]$ and $\pDA[\ge c]$
for
$\DO$, $\pDO[\le c]$ and $\pDO[\ge c]$, respectively.

\Rem\label{rem:A}
\bnum
\item A similar construction is possible for a complex manifold $X$ and
coherent $\sho_X$-modules.

\item
For any $c\in\R$, we have
\eqn
&&\ba{rcccll}
\DO[\le c]&\subset&\pDO[\le c]&\subset&\DO[\le c+\dim X/2]&\quad\text{and}
\\
\DO[\ge c+\dim X/2]&\subset&\pDO[\ge c]&\subset&\DO[\ge c].
\ea
\eneqn
\item
If $\shf$ is a Cohen-Macaulay $\sho_X$-module with $\codim(\shf)=r$,
then we have $\shf\in \pDO[-r/2]$.
\item
Assume that $A$ is a Noetherian regular integral domain of dimension $1$,
and $K$ the fraction field of $A$.
Let $\shc=\Mod_\coh(A)$.
We take as $\Tor\subset\shc$ the subcategory of torsion $A$-modules,
and as $\TF$ the subcategory of torsion free $A$-modules.
Then the t-structure 
$\bl(\pDc^{\le c}(\shc))_{c\in\R},(\pDc^{\ge c}(\shc))_{c\in\R}\br$
associated with the torsion pair
$(\Tor,\TF)$  (see \S\;\ref{sec:3})
 coincides with  the t-structure
$\bl(\pDA[\le c])_{c\in\R},(\pDA[\ge c]))_{c\in\R}\br$. 
Hence we have 
\eq&&\left\{
\ba{rcl}
\pDA[{\le n}]&=&\DA[{\le n}],\\[1ex]
\pDA[{\le n-1/2}]&=&\set{X\in \DA[{\le n}]}{K\tens_A X\in\rD^{\le n-1}(K)},\\[1ex]
\pDA[{\ge n-1/2}]&=&\DA[{\ge n}],\\[1ex]
\pDA[{\ge n}]&=&\set{X\in\DA[{\ge n}]}{\text{$H^n(X)$ is torsion free}}.
\ea\right.
\eneq
for any $n\in\Z$.

Let $\cF$ be
the quasi-abelian category of finitely generated torsion free $A$-modules.
Then
$\Db(\cF)\simeq \DA$, and the t-structure
$\bl(\pDA[\le c])_{c\in\R},(\pDA[\ge c])_{c\in\R}\br$ coincides with 
the standard t-structure of $\Db(\cF)$.
\ee
\enrem

\section{Self-dual t-structure: real case}
\label{sec:5}

\subsection{Topological perversity}
Let $X$ be a subanalytic space (cf.\ \cite[{Exercise IX.2}]{KS90}).
A subanalytic space is called smooth if it is is locally isomorphic to
a real analytic manifold as a subanalytic space.

A subanalytic stratification $X=\ssqcup[{\al\in I}]X_\al$ of $X$ is a locally finite
family of locally closed  smooth subanalytic subsets 
$\{X_\al\}_{\al\in I}$ (called strata)
such that the closure  $\ol {X_\al}$ is a union of strata for any $\al$.
A subanalytic stratification $X=\ssqcup[{\al\in I}]X_\al$ is called {\em good} 
if it satisfies the following condition:
\eq&&\parbox{70ex}{for any
$K\in\Db(\Z_X)$ such that $K\vert_{X_\al}$ has locally constant cohomologies 
for all $\al$,
$\bl\rsect_{X_\al}K\br\vert_{X_\al}$ has locally constant cohomologies 
for all $\al$.}\label{dond:good}
\eneq

Let $X=\ssqcup[{\al\in I}]X_\al$  and
$X=\ssqcup[{\al\in I'}]X_\beta'$ be two stratifications.
We say that  $X=\ssqcup[{\al\in I}]X_\al$ is finer than 
$X=\ssqcup[{\beta\in I'}]X_\beta'$ if any $X_\al$ is contained in some $X'_\beta$.
The following fact guarantees that there exist enough good stratifications:
\eq&&\parbox{70ex}{
For any locally finite family $\{Z_j\}_j$ of
locally closed subsets, there exists a good stratification
such that any $Z_j$ is a union of strata.
}
\eneq

A regular subanalytic filtration of $X$ is an increasing sequence
$$\emptyset=X_{-1}\subset\cdots\subset X_N=X$$
of closed subanalytic subsets $X_k$ of $X$ such that
$\dX_k\seteq X_k\setminus X_{k-1}$ is smooth of dimension $k$.
We say that it is a good filtration if
$\{\dX_k\}$ satisfies \eqref{dond:good}. Note that any subanalytic
stratification 
$X=\ssqcup[{\al\in I}]X_\al$ gives a regular subanalytic filtration defined by
 $X_k\seteq\ssqcupl_{\dim X_\al\le k}X_\al$.

\medskip
Let $A$ be a Noetherian regular ring. 
Let us denote by $\Mod_\Rcons(A_X)$
the category of $\R$-constructible $A_X$-modules,
and by $\DXA$ the bounded derived category of $\R$-constructible $A_X$-modules.
We denote by
$\bl(\DXA[\le c])_{c\in\R},(\DXA[\ge c])_{c\in\R}\br$ the standard t-structure of
$\DXA$, that is
\eqn
&&\DXA[\le c]=\set{K\in\DXA}{\text{$H^i(K)=0$ for any $i>c$}},\\
&&\DXA[\ge c]=\set{K\in\DXA}{\text{$H^i(K)=0$ for any $i<c$}}.
\eneqn

We define
\eq&&\hs{3ex}\left\{\hs{-.5ex}
\bns\ba{l}
\KS[\le c]=\{K\in\DXA\mid\parbox[t]{40ex}%
{$\dim\Supp(H^i(K))\le -2(i-c)$ for any $i$\},}
\\[1ex]
\KS[\ge c]=\{K\in\DXA\mid\parbox[t]{42ex}{$H^i\rsect_Z(K)=0$ for
any closed \\subanalytic subset $Z$ and
$i<c-\dim Z/2$\}.}
\ea\ens\right.
\eneq
Then we have
\Prop
The pair $\bl(\KS[\le c])_{c\in\R},(\KS[\ge c])_{c\in\R}\br$
is a t-structure on $\DXA$.
\enprop
\Proof
Indeed,
$\bl\KS[<c+1],\KS[\ge c]\br$ coincides with the t-structure
associated with the perversity
$p(n)=\lceil{c-n/2}\rceil$ (see e.g.\ \cite[Definition 10.2.1]{KS90}).
\QED

\Lemma[{\cite[Proposition 10.2.4]{KS90}}]\label{lem:KSp}
Let $K\in\DXA$ and 
let $X=\ssqcup[\al] X_\al$ be a subanalytic stratification of $X$ such that
$(\tD K)\vert_{X_\al}$ has locally constant cohomologies for any $\al$.
Then $K\in\KS[\ge c]$ if and only if
$(\rsect_{X_\al}K)_x\in\DA[\ge c-\dim X_\al/2]$ for any $\al$ and $x\in X_\al$.
\enlemma

\subsection{Self-dual t-structure : $\R$-constructible case}
As in the preceding subsection,
 $X$ is a subanalytic space and $A$ is a Noetherian regular ring.
Let $\tD$ be the duality functor
$$\tD(K)=\rhom[{A}](K,\omega_X)\quad\text{for $K\in\DXA$,}$$
where $\omega_X=\epb{a_X}A_\pt$ with 
the canonical projection $a_X\cl X\to \pt$.

For $F\in \Mod_\Rcons(A_X)$, we set
\eq
&&\ba{rcl}
\pdim(F)
&=&\sup_{m\ge0}\Bigl(\dim\set{x\in X}{\codim(F_x)=m}-m\Bigr),
\ea
\eneq
where $\codim(F_x)$ denotes the codimension of $\Supp(F_x)\subset\Spec(A)$.
Hence if $X=\ssqcup[{\al}]X_\al$ is a subanalytic stratification with connected strata and
$F\vert_{X_\al}$ is locally constant for any $\al$,
then we have
$$ \pdim(F)=\sup\set{\dim X_\al-\codim (F_{x_\al})}{F\vert_{X_\al}\not=0},$$
where $x_\al$ is a point of $X_\al$.
We understand $\pdim 0=-\infty$.

We set
\eq
\hs{4ex}
&&\ba{rl}
\pDXA[\le c]&=\set{K\in\DXA}{\text{$\pdim(H^i(K))\le -2(i-c)$ for any $i$}},
\\[1ex]
\pDXA[\ge c]&=\set{K\in\DXA}{\tD K\in\pDXA[\le -c]}.
\ea
\eneq
Note that, when $A$ is a field, they coincide with $\KS[\le c]$
and $\KS[\ge c]$.

\Lemma
Let $K\in\DXA$ and $c\in\R$.
Let $X=\ssqcup[\al] X_\al$ be a subanalytic stratification such that $K\vert_{X_\al}$ 
has locally constant cohomologies.
Then the following conditions are equivalent;
\bna
\item
$K\in \pDXA[\le c]$,
\item 
$\dim\set{x\in X}{K_x\not\in\pDA[\le c-k/2]}<k$ for any $k\in\Z$,
\item
$K_x\in\pDA[\le c-\dim X_\al/2]$ for any $\al$ and $x\in X_\al$.
\ee
\enlemma
\Proof
(a)$\Leftrightarrow$(c)\quad It is obvious that
 $K\in \pDXA[\le c]$ if and only if
$$\text{$\dim X_\al-\codim \Supp (H^i(K)_x)\le -2(i-c)$
for any $\al$, $x\in X_\al$ and  $i\in\Z$.}$$
The last condition is equivalent to
$$\codim \Supp (H^i(K_x))\ge 2(i-c+\dim X_\al/2),$$
or equivalently
$K_x\in\pDA[\le c-\dim X_\al/2]$.

\medskip\noi
(b)$\Leftrightarrow$(c)\quad (b) is equivalent to
the condition:
$$\text{for any $x\in X_\al$, $K_x\notin \pDA[\le c-k/2]$ implies $\dim X_\al<k$,}
$$
which is equivalent to the condition:
$$\text{for any $x\in X_\al$,  $\dim X_\al\ge k$
implies $K_x\in\pDA[\le c-k/2]$.}$$
It is obviously equivalent to (c).
\QED

\Lemma\label{lem crtgr}
Let $K\in\DXA$ and $c\in\R$.
Let $X=\ssqcup[\al] X_\al$ be a subanalytic stratification such that 
$(\tD K)\vert_{X_\al}$ has locally constant cohomologies.
Then the following conditions are equivalent:
\bna
\item $K\in \pDXA[\ge c]$,
\item
for any $c'\in \R$ and $M\in\pDA[\le c']$, we have
$$\rhom[A](M_X, K)\in \KS[\ge c-c'],$$
\item
$\rsect_Z(K)_x\in\pDA[{\ge c-\dim Z/2}]$ for any closed subanalytic set $Z$ and $x\in Z$,
\item
$ \bl\rsect_{X_\al}K\br_x\in \pDA[{\ge c-\dim X_\al/2}]$ for any $\al$ and $x\in X_\al$,
\item 
$\dim\set{x\in X}{\bl\rsect_{\{x\}}K\br_x\not\in \pDA[\ge c+k/2]}<k$
for any $k\in\Z_{\ge0}$.
\ee
\enlemma
\Proof
Let $i_\al\cl X_\al\to X$ be the inclusion.

\medskip\noi
(a)$\Leftrightarrow$(d)
By (a)$\Leftrightarrow$(c) in the preceding lemma,
condition (a) is equivalent to the condition:
$$\text{$(\tD K)_x\in\pDA[\le -c-\dim X_\al/2]$ for any $\al$ and $x\in X_\al$.}
$$
On the other hand, we have $\opb{i_\al}\tD K\simeq\tD[X_\al]\epb{i_\al}K$.
Hence $\epb{i_\al}K$ has locally constant cohomologies.
Since 
$$(\tD K)_x\simeq(\tD[X_\al]\epb{i_\al}K)_x\simeq\RHom[A]\bl(\epb{i_\al}K)_x,A\br[\dim X_\al],$$
the above condition is equivalent to
$$\RHom[A]\bl(\epb{i_\al}K)_x,A\br\in\pDA[{\le -c+\dim X_\al/2}],$$
which is again equivalent to
$(\epb{i_\al}K)_x\in\pDA[{\ge c-\dim X_\al/2}]$.

\medskip\noi
(a)$\Leftrightarrow$(e)\quad (a) is equivalent to
$\tD K\in\pDX[\le -c]$.
By the preceding lemma, it is equivalent to the condition:
$\dim\set{x\in X}{(\tD K)_x\not\in\pDA[\le -c-k/2]}<k$ for any $k\in\Z_{\ge0}$.
Since $(\tD K)_x\simeq\tD[A]\bl(\rsect_{\{x\}}K)_x\br$, the condition
$(\tD K)_x\not\in\pDA[\le -c-k/2]$
is equivalent to
$(\rsect_{\{x\}}K)_x\not\in\pDA[\ge c+k/2]$.

\medskip\noi
(d)$\Leftrightarrow$(b)\ 
Condition (d) is equivalent to
\eq
&&\parbox{70ex}{$\RHom\bl M,(\rsect_{X_\al}K)_x\br
\in \DA[{\ge c-\dim X_\al/2-c'}]$ 
for any $M\in\pDA[\le c']$, $\al$ and $x\in X_\al$.}\label{cond:int}
\eneq
Since $\RHom\bl M,(\rsect_{X_\al}K)_x\br\simeq
\bl\rsect_{X_\al}\rhom[A](M_X,K)\br_x$,
the last condition \eqref{cond:int} is equivalent to
(b) by Lemma~\ref{lem:KSp}.

\medskip\noi
(c)$\Rightarrow$(d) is obvious.

\medskip\noi
(b)$\Rightarrow$(c)\ 
For any $c'\in \R$ and $M\in\pDA[\le c']$, we have
$\bl\rsect_Z\rhom[A](M_X,K)\br_x\in\DA[\ge c-c'-\dim Z/2]$.
Since 
$\RHom[A]\bl M,(\rsect_ZK)_x\br\simeq \bl\rsect_Z\rhom[A](M_X,K)\br_x$,
we have  (c)

\QED

We shall prove the following theorem in several steps.
\Th\label{th:main}
$\bl (\pDXA[\le c])_{c\in\R},(\pDXA[\ge c])_{c\in\R}\br$ is a t-structure on \/
$\DXA$.
\enth

It is obvious that it satisfies conditions \eqref{ax:1} and \eqref{ax:2}
in Definition~\ref{axiom}.
Let us show \eqref{ax:3}.

\Lemma
For $c\in\R$ and $K\in\pDXA[\le c]$ and $L\in\pDXA[\ge c']$, we have
$$\rhom(K,L)\in \Der^{\ge c'-c}_\Rcons(A_X).$$
\enlemma
\Proof
Let us take a \good regular subanalytic filtration of $X$
$$\emptyset=X_{-1}\subset\cdots\subset X_N=X$$
such that $K$ and $L$ have locally constant cohomologies on
each $\dX_k\seteq X_k\setminus X_{k-1}$.
We may assume that $\dX_k$ is smooth of dimension $k$.

Let $i_k\cl \dX_k\to X$ 
be the inclusion.

Let us first show that
\begin{align}
\text{$\epb{i_k}\rhom(K,L)\simeq 
\rhom(\opb{i_k}K,\epb{i_k}L)$
belongs to $\Der^{\ge c'-c}_\Rcons(A_{\dX_k})$.}
\label{eq:smoothpt}
\end{align}

Since $\opb{i_k}K$, $\epb{i_k}L$
have locally constant cohomologies, 
$$\bl\epb{i_k}\rhom(K,L)\br_x\simeq
\rHom[A]\bl(\opb{i_k}K)_x,(\epb{i_k}L)_x\br$$
for any $x\in \dX_k$.
Hence it is enough to show that
\eq
&&\rHom[A]\bl(\opb{i_k}K)_x,(\epb{i_k}L)_x\br\in\Der^{\ge c'-c}_\Rcons(A).
\label{eq:pointvan}
\eneq
It follows from Corollary~\ref{cor:hom} and
$(\opb{i_k}K)_x\in \pDA[\le c-k/2]$ and $(\epb{i_k}L)_x\in \pDA[\ge c'-k/2]$.

Now we shall show that
$$\rsect_{X_k}\rhom(K,L)\in \Der^{\ge c'-c}_\Rcons(A_X)$$
by induction on $k$.
By the induction hypothesis
$\rsect_{X_{k-1}}\rhom(K,L)\in \Der^{\ge c'-c}_\Rcons(A_X)$.
We have the distinguished triangle
$$\rsect_{X_{k-1}}\rhom(K,L)\To
\rsect_{X_k}\rhom(K,L)\To\rsect_{\dX_k}\rhom(K,L)\tone$$
Since
$\rsect_{\dX_k}\rhom(K,L)\simeq\RR i_{k\;*}\epb{i_k}\rhom(K,L)$
belongs to $\Der^{\ge c'-c}_\Rcons(A_X)$,
we obtain $\rsect_{X_k}\rhom(K,L)\in \Der^{\ge c'-c}_\Rcons(A_X)$.
\QED

Now we shall show Definition~\ref{axiom} \eqref{ax:4} in a special case.
\Lemma\label{lem:smooth}
Let us assume that $X$ is a smooth subanalytic space, and $c\in\R$.
Let $K\in\DXA$ and assume that $K$ has locally constant cohomologies.
Then there exists a \dt
$$K'\To K\To K''\tone$$
with $K'\in\pDXA[\le c]$ and $K''\in\pDXA[>c]$.
Moreover $K'$ and $K''$ have locally constant cohomologies.
\enlemma
\Proof
Let us show it in three steps.

\noi
(i) There exists locally such a \dt.

\smallskip
Indeed, for any $x\in X$, there exist an open neighborhood $U$ of $x$ and
 $M\in\DA$ such that $K\vert_U\simeq M_U$.
Take a \dt
$M'\to M\to M''\tone$
such that $M'\in\pDA[\le c-\dim X/2]$ and
$M''\in \pDA[>c-\dim X/2]$.
Then $M'_{U}\to M_U\to M''_U\tone$ gives a desired \dt.

\medskip\noi
(ii) If $U_i$ is an open subset of
$X$ and $K'_i\To K\vert_{U_i}\To K''_i\tone$ is a \dt with
$K'_i\in\pD^{\le c}_\Rcons(A_{U_i})$ and $K''_i\in\pD^{>c}_\Rcons(A_{U_i})$  ($i=1,2$),
then there exists a \dt
$K'\To K\vert_{U_1\cup U_2}\To K''\tone$ with 
$K'\in\pD^{\le c}_\Rcons(A_{U_1\cup U_2})$ and $K''\in\pD^{>c}_\Rcons(A_{U_1\cup U_2})$.

\smallskip
Indeed, by the uniqueness of such a \dt, we have
$K'_1\vert_{U_1\cap U_2}\simeq K'_2\vert_{U_1\cap U_2}$.
Denote it by $K_0\in \Db(A_{U_1\cap U_2})$.
Let $i_0\cl U_1\cap U_2\to U_1\cup U_2$
and $i_k\cl U_k\to  U_1\cup U_2$($k=1,2$) be the open inclusions.
Then embed
a morphism
$\eim{i_0}K_0\to \eim {i_1} K'_1\soplus\eim{i_2} K'_2$ into a distinguished triangle
$$\eim{i_0}K_0\to \eim {i_1} K'_1\soplus\eim{i_2} K'_2\to K'\tone.$$
Then $K'\vert_{U_k}\simeq K'_k$.
Since the composition
$\eim{i_0}K_0\to \eim {i_1} K'_1\soplus\eim{i_2} K'_2\to K\vert_{U_1\cup U_2}$ vanishes, the morphism $ \eim {i_1} K'_1\soplus\eim{i_2} K'_2\to K\vert_{U_1\cup U_2}$
factors through $K'$.
Hence, there exists a morphism $K'\to K\vert_{U_1\cup U_2}$ which extends
$K'_i\to K\vert_{U_i}$ ($i=1,2$).
Then, embedding this morphism into a \dt
$K'\To K\vert_{U_1\cup U_2}\To K''\tone$,
we obtain the desired \dt.

\medskip\noi
(iii)\ By (i) and (ii), there exist an increasing sequence of open subsets 
$\{U_n\}_{n\in\Z_{\ge0}}$ with $X=\scup[{n\in\Z_{\ge0}}] U_n$ and a \dt
$K_n\to K\vert_{U_n}\to K''_n\tone$ with
$K'_n\in\pD^{\le c}_\Rcons(A_{U_n})$ and $K''_n\in\pD^{>c}_\Rcons(A_{U_n})$.
Let $i_n\cl U_n\to X$ be the inclusion.
By the uniqueness of such \dts,
we have
$K'_{n+1}\vert_{U_n}\simeq K'_n$. Hence,
we have $\beta_n\cl\eim{i_n}K'_n\to\eim{i_{n+1}}K'_{n+1}$.
Let $K'$ be the hocolim of the inductive system
$\{\eim{i_n}K'_n\}_{n\in\Z_{\ge0}}$, that is,
the third term of a \dt
$$\soplus_{n\in\Z_{\ge0}}\eim{i_n}K'_n\To[\ f\ ] \soplus_{n\in\Z_{\ge0}}\eim{i_n}K'_n\To K'\tone.$$
Here $f$ is given such that the following diagram commutes for any $a\in\Z_{`\ge0}$:

$$\xymatrix@C=15ex@R=3ex{
\eim{i_a}K'_a\ar[r]^-{\id_{\eim{i_a}K'_a}\oplus\, (-\beta_{a})}\ar[d]&
\eim{i_a}K'_a\oplus\eim{i_{a+1}}K'_{a+1}\ar[d]\\
\soplus_{n\in\Z_{\ge0}}\eim{i_n}K'_n\ar[r]^f
&\soplus_{n\in\Z_{\ge0}}\eim{i_n}K'_n.
}
$$

Then $K'\vert_{U_n}\simeq K'_n$.
Since the composition
$$\soplus_{n\in\Z_{\ge0}}\eim{i_n}K'_n\To[\ f\ ] \soplus_{n\in\Z_{\ge0}}\eim{i_n}K'_n\To 
K$$ vanishes, the morphism 
$\soplus_{n\in\Z_{\ge0}}\eim{i_n}K'_n\To K$ factors through $K'$.
Hence there is a morphism $K'\to K$ which extends
$\eim{i_n}K'_n\to K$.
Hence, embedding this morphism into a \dt
$$K'\To K\To K''\tone,$$
we obtain the desired \dt.
\QED
Finally we shall complete the proof of Definition~\ref{axiom} \eqref{ax:4}.
\Lemma
Let $K\in\DXA$ and $c\in\R$. Then there exists a distinguished triangle
$$K'\To K\To K''\tone$$
with $K'\in\pDXA[\le c]$ and $K''\in\pDXA[>c]$.
\enlemma
\Proof
Let us take a \good regular subanalytic filtration of $X$
$$\emptyset=X_{-1}\subset\cdots\subset X_N=X$$
such that $K$ has locally constant cohomologies on
each $\dX_k\seteq X_k\setminus X_{k-1}$.
We may assume that $\dX_k$ is a smooth subanalytic space of dimension $k$.
We shall prove that
\begin{equation*}
\addtocounter{equation}{1}
\setcounter{mycounter}{\value{section}}
\setcounter{mycounterf}{\value{subsection}}
\setcounter{mycounters}{\value{equation}}
\hs{-3ex}(\themycounter.\themycounterf.\themycounters)_{k}\hs{1ex}
\left\{\parbox{70ex}{there exists a distinguished triangle\\
\hs{10ex}$K'\To K\vert_{X\setminus X_{k}}\To K''\tone$\\
with  $K\in\pD^{\le c}(A_{X\setminus X_k})$ and $K''\in\pDXA[>c]$.
Moreover, $K'\vert_{\dX_j}$ and $K''\vert_{\dX_j}$ have
locally constant cohomologies for $j>k$.}
\right.
\end{equation*}
by the descending induction on $k$.

Assuming $(\themycounter.\themycounterf.\themycounters)_{k}$,
we shall show $(\themycounter.\themycounterf.\themycounters)_{k-1}$.
Let $K'\To K\vert_{X\setminus X_{k}}\To K''\tone$ be a \dt
as in $(\themycounter.\themycounterf.\themycounters)_{k}$. Let 
$j\cl X\setminus X_k\to X\setminus X_{k-1}$ be the open embedding
and $i\cl \dX_k\to X\setminus X_{k-1}$ the closed embedding.
The morphism $K'\to K\vert_{X\setminus X_{k}}$ induces
$\eim{j} K'\to  K\vert_{X\setminus X_{k-1}}$.
We embed it into a \dt  in $\Der^{\mathrm{b}}_\Rcons(A_{X\setminus X_{k-1}})$
$$\eim{j} K'\to  K\vert_{X\setminus X_{k-1}}\to L\tone.$$
By Lemma~\ref{lem:smooth}, there exists a \dt
\eq
&&L'\To \epb{i}L\To L''\To[\;+1\;]
\label{eq:L}
\eneq
with $L'\in\pD^{\le c}_\Rcons(A_{\dX_k})$  and $L''\in\pD^{> c}_\Rcons(A_{\dX_k})$.
We embed the composition $\eim{i}L'\to\eim{i}\epb{i}L\to L$ into a \dt
\eq
&&\eim{i}L'\to L\to \tK''\tone
\label{eq:tK''}
\eneq
Finally we embed the composition $ K\vert_{X\setminus X_{k-1}}\to L\to \tK''$
into a \dt
$$\tK'\to  K\vert_{X\setminus X_{k-1}}\to \tK''\tone.$$
Let us show that
$$\text%
{$\tK'\in \pD^{\le c}_\Rcons(A_{X\setminus X_{k-1}})$ and 
$\tK''\in \pD^{> c}_\Rcons(A_{X\setminus X_{k-1}})$.}
$$
By the construction, we have
$\tK''\vert_{X\setminus X_{k}}\simeq L\vert_{X\setminus X_{k}}\simeq K''$
and $\tK'\vert_{X\setminus X_{k}}\simeq K'$.
Hence it is enough to show that
$\opb{i}\tK'\in \pD^{\le c}_\Rcons(A_{\dX_k})$ and
$\epb{i}\tK''\in \pD^{>c}_\Rcons(A_{\dX_k})$.
Applying the functor $\epb{i}$ to \eqref{eq:tK''}, we obtain
a \dt
$$L'\To \epb{i}L\To\epb{i}\tK''\tone.$$
By the \dt \eqref{eq:L}, we have
$\epb{i}\tK''\simeq L''\in \pD^{>c}_\Rcons(A_{\dX_k})$.

\medskip
By the octahedral axiom for triangulated category,
 we have a diagram
\eqn
\xymatrix@C=6ex{&\tK'\ar@{.>}[dr]\ar[ddl]\\
{\ \eim{j}K'\ }\ar@{.>}[ur]\ar[d]&&{\ \eim{i}L'\ }\ar[ll]_{+1}\ar[ddl]\\
K\vert_{X\setminus X_{k-1}}\ar[dr]\ar[rr]&&{\ \tK''\ }\ar[u]_{+1}\ar[uul]|(.65){+1}\\
&L\ar[ur]\ar[uul]|(.63){+1}
}
\eneqn
and a distinguished triangle 
$$\eim{j}K'\To\tK'\To\eim{i}L'\tone.$$
It implies
$\opb{i}\tK'\simeq L'\in \pD^{\le c}_\Rcons(A_{\dX_k})$.
\QED

This completes the proof of Theorem~\ref{th:main}.

\bigskip
Recall the full subcategory of $\DXA$:
$$\pDXA[{[a,b]}]\seteq\pDXA[\le b]\cap\pDXA[\ge a]$$
for $a\le b$.

\Prop
Assume that $a, b\in \R$ satisfy $a\le b<a+1$. Then
$X\supset U\longmapsto \pD^{[a,b]}_\Rcons((A_U)$ is a stack on $X$.
\enprop
\Proof
(i)\ Let $K,L\in\pDXA[{[a,b]}]$. Since $\rhom[A](K,L)\in\DXA[\ge a-b]
=\DXA[\ge0]$,
the presheaf
$U\mapsto\Hom_{\pD^{[a,b]}_\Rcons(A_U)}(K\vert_U, L\vert_U)\simeq 
\sect\bl U;H^0(\rhom[A](K,L))\br$ is a sheaf.
Hence, $U\mapsto \pD^{[a,b]}_\Rcons(A_U)$ is a separated prestack on $X$.

\medskip\noi
(ii)\ Let us show the following statement:
{\addtolength{\mylength}{-4ex}
\eqn
&&\parbox{\mylength}%
{Let $U_1$ and $U_2$ be open
subsets of $X$ such that $X=U_1\cup U_2$,
and $K_k\in \pD^{[a,b]}_\Rcons(A_{U_k})$ ($k=1,2$). Assume that
$K_1\vert_{U_1\cap U_2}\simeq K_2\vert_{U_1\cap U_2}$.
Then there exists 
$K\in \pDXA[{[a,b]}]$ such that $K\vert_{U_k}\simeq K_k$ ($k=1,2$).
}
\eneqn
}
\smallskip
Set $U_0=U_1\cap U_2$ and $K_0=K_1\vert_{U_1\cap U_2}\simeq K_2\vert_{U_1\cap U_2}
\in\pD^{[a,b]}_\Rcons(A_{U_0})$.
Let $j_k\cl U_k\to X$ be the open inclusion ($k=0,1,2$).
Then we have
$\beta_k\cl \eim{j_0}(K_0) \to \eim{j_k}K_k$ ($k=1,2$).
We embed the morphism $(\beta_1,\beta_2)\cl
\eim{j_0}(K_0) \to \eim{j_1}K_1\soplus\eim{j_2}K_2$ into a \dt
$$\eim{j_0}(K_0) \To \eim{j_1}K_1\soplus\eim{j_2}K_2\To K\tone.$$
Then $K$ satisfies the desired condition.

\medskip\noi
(iii)\ Let us show the following statement:
{\addtolength{\mylength}{-4ex}
\eqn
&&\parbox{\mylength}%
{Let $\{U_n\}_{n\in\Z_{\ge0}}$ be an increasing
sequence of open subsets of $X$ such that 
$X=\scup[n\in\Z_{\ge0}]U_n$.
Let $K_n\in \pD^{[a,b]}_\Rcons(A_{U_n})$ ($n\in\Z_{\ge0}$) and
$K_{n+1}\vert_{U_n}\simeq K_n$.
Then there exists 
$K\in \pDXA[{[a,b]}]$ such that $K\vert_{U_n}\simeq K_n$ ($n\in\Z_{\ge0}$).
}
\eneqn
}
\smallskip
The proof is similar to 
proof of Lemma~\ref{lem:smooth}.
Let $j_n\cl U_n\to X$ be the
open inclusion, and let $\eim{j_n}K_n\to\eim{j_{n+1}}K_{n+1}$
be the morphism induced by the isomorphism $K_{n+1}\vert_{U_n}\simeq K_n$.
Let $K$ be a hocolim of the inductive system $\{\eim{j_n}K_n\}_{n\in\Z_{\ge0}}$.
Then $K\in\pDXA[{[a,b]}]$ satisfies the desired condition.

(iv) By (i)--(iii), we conclude that $U\mapsto \pD^{[a,b]}_\Rcons(A_U)$ is a stack on $X$.
\QED

\Prop\label{prop:funct}
Let $f\cl X\to Y$ be a morphism of subanalytic spaces, and $d\in \Z_{\ge0}$.
Assume that $\dim f^{-1}(y)\le d$ for any $y\in Y$.
Then
\bnum
\item
If $G\in\pDYA[\le c]$, then $\opb{f}G\in\pDXA[\le c+d/2]$.
\item
If $G\in\pDYA[\ge c]$, then $\epb{f}G\in\pDXA[\ge c-d/2]$.
\item
If $F\in\pDXA[\ge c]$ and $\roim{f}F\in\DYA$, then 
$\roim{f}F\in\pDYA[\ge c-d/2]$.
\item
If $F\in\pDXA[\le c]$ and $\reim{f}F\in\DYA$, then 
$\reim{f}F\in\pDYA[\le c+d/2]$.
\ee
\enprop
\Proof
(i)\  Assume $G\in\pDYA[\le c]$.
Then 
\eqn
&&\dim\set{x\in X}{(\opb{f}G)_x\notin\pDA[\le c+d/2-k/2]}\\
&&\hs{10ex}=\dim f^{-1}\left(\set{y\in Y}{G_y\notin\pDA[\le c+d/2-k/2]}\right)\\
&&\hs{10ex}\le\dim\set{y\in Y}{G_y\notin\pDA[\le c+d/2-k/2]}+d\\
&&\hs{10ex}<(k-d)+d=k.
\eneqn

\medskip
\noi 
(ii) follows from (i) by the duality.

\medskip
\noi
(iii)\ For any $G\in\pDYA[<c-d/2]$, then
$$\Hom_{\DYA}(G,\roim{f}F)\simeq\Hom_{\DXA}(\opb{f}G,F)$$
vanishes because $\opb{f}G\in\pDXA[<c]$ by (i).
Hence $\roim{f}F\in\pDYA[\ge c-d/2]$ by \eqref{eq:orth}.

Similarly, (iv) follows from (ii).
\QED

We shall give relations between the two t-structures:\\
 $\bl(\KS[\le c])_{c\in\R},(\KS[\ge c])_{c\in\R}\br$ and
$\bl (\pDXA[\le c])_{c\in\R},(\pDXA[\ge c])_{c\in\R}\br$.

\Lemma\label{lem:KSP}
Let $K\in\DXA$ and $c\in\R$.
\bnum
\item
The following conditions are equivalent:
\bna
\item
$K\in\pDXA[\le c]$,
\item
for any $c'\in \R$ and $M\in\pDA[\ge c']$, we have
$$\rhom[A](K,M\tens\omega_X)\in \KS[\ge c'-c].$$
\ee
\item
The following conditions are equivalent:
\bna
\item
$K\in\pDXA[\ge c]$,
\item
for any $c'\in \R$ and $M\in\pDA[\le c']$, we have
$$\rhom[A](M_X, K)\in \KS[\ge c-c'].$$
\ee
\ee
\enlemma
\Proof (ii) is already proved in Lemma~\ref{lem crtgr}.
(i) follows from (ii) because
\eqn\rhom[A](K,M\tens\omega_X)
&\simeq&\rhom[A]\bl(\tD(M\tens\omega_X),\tD K\br\\
&\simeq&\rhom[A]\bl(\tD[A]M)_X,\tD K\br,\eneqn
where $\tD[A]M\seteq\rHom[A](M,A)$.
\QED
\Lemma\label{lem:extpr}
Let $X$ and $Y$ be subanalytic spaces.
Let $K\in \pDXA[\ge c]$ and $L\in\pDYA[\ge c']$.
Then we have
$K\letens L\in\ks_\Rcons^{\ge c+c'}(A_{X\times Y})$.
\enlemma
\Proof
Let $X=\ssqcup[\al] X_\al$ and $Y=\ssqcup[\beta] Y_\beta$
be good subanalytic stratification such that $K\vert_{X_\al}$ and $L\vert_{Y_\beta}$
are locally constant.
Then we have
$(\rsect_{X_\al}K)_x\in\pDA[\ge c-\dim X_\al/2]$
and $(\rsect_{Y_\beta}L)_y\in\pDA[\ge c'-\dim Y_\beta/2]$
for $x\in X_\al$ and $y\in Y_\beta$.
Hence by Proposition~\ref{prop:inner} \eqref{item:3},
we have $$\bl\rsect_{X_\al\times Y_\beta}(K\letens L)\br_{(x,y)}
\simeq(\rsect_{X_\al}K)_x\ltens (\rsect_{Y_\beta}L)_y\in \DA[\ge c+c'-\dim(X_\al\times Y_\beta)/2].$$ 
Hence we obtain
$K\letens L\in\ks_\Rcons^{\ge c+c'}(A_{X\times Y})$.
\QED

\Rem
We have
\eqn
\KS[\le c]&\subset& \pDXA[\le c],\\[1ex]
\pDXA[\ge c]&\subset& \KS[\ge c].
\eneqn
\enrem

\section{Self-dual t-structure: complex analytic variety  case}
\label{sec:6}
\subsection{Middle perversity in the complex case}
Let $X$ be a complex analytic space.
We denote by $\dc X$ the dimension of $X$.
Hence we have
 $\dc X=(\dim X_\R)/2$ where $X_\R$ is the underlying subanalytic space.
For a complex submanifold $Y$ of a complex manifold $X$, we denote by $\codc Y$ the codimension 
of $Y$ as complex manifolds.
We sometimes write $d_X$ for $\dc X$.

Let $\DXAC$ be the bounded derived category of 
the abelian category of sheaves of $A$-modules with
$\C$-constructible cohomologies.
Then it is a full subcategory of $\DXA$ and it is easy to see that
the self-dual t-structure on $\DXA$ induces a self-dual
t-structures on $\DXAC$.
More precisely, if we define
\eqn
&&\pDXAC[\le c]\seteq\DXAC\cap \pDXA[\le c]\quad\text{and}\\
&&\pDXAC[\ge c]\seteq\DXAC\cap \pDXA[\ge c],
\eneqn
then $\bl (\pDXAC[\le c])_{c\in\C},(\pDXAC[\ge c])_{c\in\C}\br$ is a t-structure on
$\DXAC$.

\noi
Similarly, the t-structure
$\bl (\KS[\le c])_{c\in\C},(\KS[\ge c])_{c\in\C}\br$ 
induces 
the t-structure 
$\bl (\KSC[\le c])_{c\in\C},(\KSC[\ge c])_{c\in\C}\br$
on $\DXAC$.

Note that the t-structure $\bl\KSC[\le 0],\KSC[\ge 0]\br$ in the original sense
is denoted by $\bl{}^{p}\mathbf{D}_{\C\text{-}c}^{\le0}(X),
{}^{p}\mathbf{D}_{\C\text{-}c}^{\ge0}(X)\br$
 in \cite[\S 10.3]{KS90}.

In \cite[\S 10.3]{KS90},
various properties of $\bl\KSC[\le 0],\KSC[\ge 0]\br$ are studied. 
By using Lemma~\ref{lem:KSP},
we have similar properties for $\bl (\pDXAC[\le c])_{c\in\C},(\pDXAC[\ge c])_{c\in\C}\br$ as we explain in the next subsection.

\subsection{Microlocal characterization} 
Let us assume that $X$ is a complex manifold.
Let $K\in\DXAC$. Then its microsupport $\Ss(K)$ is a Lagrangian complex analytic subset of
the cotangent bundle $T^*X$ (see \cite{KS90}).

A point $p$ of $\Ss(K)$ is called {\em good} if $\Ss(K)$ is equal to
the conormal bundle $T^*_YX$ on a neighborhood of $p$
for some locally closed complex submanifold $Y$ of $X$.
The complement of the set of good points of $\Ss(K)$
is a nowhere dense closed complex analytic subset of $\Ss(K)$.
For such a good point $p$ of $\Ss(K)$,
there exists $L\in\DA$ such that
$K$ is micro-locally isomorphic to $L_Y[-\codc Y]$ on a neighborhood of $p$.
We call $L$ the type of $K$ at $p$.
(Note that $L$ is called the type of $K$ 
at $p$ with shift $0$ in \cite[\S 10.3]{KS90}.)

The type can be calculated by the vanishing cycle functor.
If we take a holomorphic function such that
$f\vert_Y=0$ and $df(x_0)=p$,
then we have
$\vphi_f(K)_{x_0}\simeq L[-\codc Y]$.
Here, $x_0\in X$ is the image of $p$ by the projection $T^*X\to X$, and
$\vphi_f$ is the vanishing cycle functor
from $\Der^{\mathrm{b}}_\Ccons(A_X)$ to $\Der^{\mathrm{b}}_\Ccons(A_{f^{-1}(0)})$.
Note that we have an isomorphism $$\vphi_f(K)\simeq
\rsect_{\{x\mid \Re(f(x))\ge0\}}(K)\vert_{f^{-1}(0)}.$$

The following theorem is proved in \cite[\S 10.3]{KS90}.
\Th[{\cite[Theorem 10.3.2]{KS90}}]
Let $K\in \DXAC$. Then the following conditions are equivalent.
\bna
\item
$K\in \KSC[\le c]$ \ro resp.\ $K\in \KSC[\ge c]$\rf,
\item
the type of $K$ at any good point of\/ $\Ss(K)$
belongs to $\DA[\le c-d_X]$ \ro resp.\ belongs to $\DA[\ge c-d_X]$\rf.
\ee
\enth

As a corollary of this theorem, we can derive the following microlocal characterization of $\bl (\pDXAC[\le c])_{c\in\C},(\pDXAC[\ge c])_{c\in\C}\br$.

\Th\label{th:micro}
Let $K\in \DXAC$. Then the following conditions are equivalent.
\bna
\item
$K\in \pDXAC[\le c]$ \ro resp.\ $K\in \pDXAC[\ge c]$\rf,
\item
the type of $K$ at any good point of\/ $\Ss(K)$
belongs to $\pDA[\le c-d_X]$ \ro resp.\ belongs to $\pDA[\ge c-d_X]$\rf.
\ee
\enth
\Proof
Let us assume that $K\in \pDXA[\ge c]$.
Then for any $M\in\pDA[\le c']$, we have
$\rhom[A](M_X,K)\in\KSC[\ge c-c']$.
Let $L$ be the type of $K$ at  a good point $p$ of $\Ss(K)$.
Then, $\rhom[A](M_X,K)$ has type $\RHom[A](M,L)$ at $p$.
Hence,  the preceding theorem implies
$\RHom[A](M,L)\in\DA[\ge c-c'-d_X]$.
Since this holds for any $M\in\pDA[\le c']$,
we conclude
$L\in\pDA[\ge c-d_X]$.
The converse can be proved similarly.

\medskip
The case of $ \pDXAC[\le c]$ can be derived from the case of
$ \pDXAC[\ge c]$ by duality.
Condition $K\in \pDXAC[\le c]$ is equivalent to
$\tD(K) \in \pDXAC[\ge -c]$.
Let $L$ be the type of $K$ at a good point $p$ of $\Ss(K)$.
Then, $\tD(K)$ has type $\tD[A](L)[2d_X]$ at $p$.
Then it is enough to notice that
$\tD[A](L)[2d_X]\in\pDA[\le -c-d_X]$,
if and only if
$L\in \pDA[\ge c-d_X]$.
\QED
The following proposition can be proved similarly.
\Prop
Let $Y$ be a closed complex submanifold
of a complex manifold $X$.
Then we have
\bnum
\item The functor
 $\nu_Y\cl\DXAC\to \Der^{\mathrm{b}}_\Ccons(A_{T_YX})$ sends 
$\pDXAC[\le c]$ to $\pD^{\le c}_\Ccons(A_{T_YX})$ and
$\pDXAC[\ge c]$ to $\pD^{\ge c}_\Ccons(A_{T_YX})$,
\item The microlocalization functor
 $\mu_Y\cl\DXAC\to \Der^{\mathrm{b}}_\Ccons(A_{T^*_YX})$ sends 
$\pDXAC[\le c]$ to $\pD^{\le c +\codc Y}_\Ccons(A_{T^*_YX})$ and
$\pDXAC[\ge c]$ to $\pD^{\ge c+\codc Y}_\Ccons(A_{T^*_YX})$.
\ee
\enprop
\Proof
Since the proof is similar, we show only (ii).
Let $K\in \pDXAC[\ge c]$. Then, for any
$M\in\pDA[\le c']$, we have
$\rhom[A](M_X,K)\in\KSC[\ge c-c']$.
Hence \cite[{Prop.~10.3.19}]{KS90} implies that
$\mu_Y\bl\rhom[A](M_{X},K)\br\in \KSXY[\ge c-c'+\codc Y]$.
Since we have
$$\rhom[A](M_{T^*_YX},\mu_YK)\simeq\mu_Y\bl\rhom[A](M,K)\br,$$
we obtain 
$\mu_YK\in \pD^{\ge c +\codc Y}_\Ccons(A_{T^*_YX})$.

Assume now that $K\in \pDXAC[\le c]$.
Then we have $\tD K\in \pDXAC[\ge -c]$.
Since \cite[{Prop.~8.4.13}]{KS90}
implies $\tD[T^*_YX](\mu_YK)\simeq(\mu_Y\tD K)^a[2\codc Y]$, we obtain
$$\tD[T^*_YX](\mu_YK)\in \pD^{\ge -c -\codc Y}_\Ccons(A_{T^*_YX}).$$
Hence $\mu_YK\in \pD^{\le c+\codc Y}_\Ccons(A_{T^*_YX})$.
\QED

The following theorem is proved in \cite[\S 10.3]{KS90}.
\Th[{\cite[Corollary 10.3.20]{KS90}}]\label{th:KS}
Let $K\in\KSC[\le c]$ and $L\in\KSC[\ge c']$.
Then $\muhom(K,L)\in\KST[\ge c'-c+d_X]$.
\enth

As its corollary we obtain the following result.
\Th Let $K\in\DXAC$ and $L\in\DXAC$.
\bnum
\item
If $K\in\pDXAC[\le c]$ and  $L\in\pDXAC[\ge c']$,
then we have $$\muhom(K,L)\in\KST[\ge c'-c+d_X].$$
\item
If $K\in\KSC[\le c]$ and  $L\in\pDXAC[\ge c']$,
then we have $$\muhom(K,L)\in\pDXT[\ge c'-c+d_X].$$
\ee
\enth
\Proof
(i)\ By Lemma~\ref{lem:extpr}, we have
$L\letens\tD K\in\KS[\ge c'-c]$.
Let $\Delta_X$ be the diagonal set of $X\times X$.
Then we have
$\muhom(K,L)=\mu_{\Delta_X}(L\letens \tD K)
\in \KSC[\ge c'-c+d_X]$ by \cite[{Proposition 10.3.19}]{KS90}.

\medskip\noi
(ii)\ For any $M\in\pDA[\le c'']$, we have 
$\rhom(M_X,L)\in \KSC[\ge c'-c'']$.
Hence
$$\rhom(M_{T^*X},\muhom(K,L))\simeq
\muhom(K,\rhom(M_X,L))$$ belongs to
$\KST[\ge c'-c''-c+d_X]$ by Theorem~\ref{th:KS}.
Hence we conclude $\muhom(K,L)\in \pDXT[\ge c'-c+d_X]$
by Lemma~\ref{lem:KSP}.

\QED

\Ex
Assume that $2$ acts injectively on $A$.
Let $M$ be a finitely generated projective $A$-module.
Let $X=\C^3$ and $S=\set{(x,y,z)\in X}{x^2+y^2+z^2=0}$.
Let $j\cl X\setminus\{0\}\to X$ be the inclusion.
Since $S\setminus\{0\}$ is homeomorphic to the product of
$\R$ and  the $3$-dimensional real
projective space $\mathbb{P}^3(\R)$,  we have
$$\bl\roim{j}\opb{j}(M_S)\br_0\simeq\rsect(S\setminus\{0\}; M_S)
\simeq M\soplus (M/2M)[-2]\soplus M[-3],$$
and $\rsect_{\{0\}}(M_S)_0\simeq (M/2M)[-3]\soplus M[-4]$.
Hence we have 
$$M_S\in\pDXAC[2],$$
and a \dt
$$M_{0}[-1]\to\reim{j}\opb{j}(M_S)\to M_S\tone.$$
Hence we obtain
\eqn
\reim{j}\opb{j}(M_S)&\in &\pDXAC[{[1,2]}],\\
\ptau^{\ge2}\reim{j}\opb{j}(M_S)&\simeq&M_S, \\
\ptau^{<2}\reim{j}\opb{j}(M_S)&\simeq& M_0[-1]\in \pDXAC[1].
\eneqn
Here $\ptau$ denotes the truncation functor of the t-structure $\pDXAC$.

By the duality, we have
\eqn
\roim{j}\opb{j}(M_S)&\in &\pDXAC[{[2,3]}],\\
\ptau^{>2}\roim{j}\opb{j}(M_S)&\simeq& M_0[-3]\in \pDXAC[3].
\eneqn
Hence we obtain a \dt
$$\ptau^{\le 2}\roim{j}\opb{j}(M_S)\to \roim{j}\opb{j}(M_S)\to M_0[-3]\tone.$$
The canonical morphism $\reim{j}\opb{j}(M_S)\to \roim{j}\opb{j}(M_S)$ decomposes as
$$\xymatrix
{\reim{j}\opb{j}(M_S)\ar[r]\ar[d]&\roim{j}\opb{j}(M_S)\\
M_S\ar[r]&\ptau^{\le 2}\roim{j}\opb{j}(M_S)\ar[u]
}
$$
and the bottom arrow is embedded into a \dt
$$M_S\to \ptau^{\le 2}\roim{j}\opb{j}(M_S)\to(M/2M)_{\{0\}}[-2]\tone.$$
Note that $(M/2M)_{\{0\}}[-2]\in\pDXAC[3/2]$.
Hence $M_S\to \ptau^{\le 2}\roim{j}\opb{j}(M_S)$
is a monomorphism and an epimorphism in the quasi-abelian category
$\pDXAC[2]$.
Moreover, we have an exact sequence
$$0\to M_S\to \ptau^{\le 2}\roim{j}\opb{j}(M_S)\to(M/2M)_{\{0\}}[-2]\to0$$
in the abelian category $\pDXAC[{[3/2,\,2]}]$ and
an exact sequence
$$0\to(M/2M)[-3]_{\{0\}}\to M_S\to \ptau^{\le 2}\roim{j}\opb{j}(M_S)\to0$$
in the abelian category $\pDXAC[{[2,\,5/2]}]$.
Note that we have an isomorphism of \dts
$$
\xymatrix@R=3ex{
\vphi_x(M_S)\ar[r]\ar[d]^\bwr&\vphi_x\bl\ptau^{\le 2}\roim{j}\opb{j}(M_S)\br
\ar[r]\ar[d]^\bwr&
\vphi_x\bl(M/2M)_{\{0\}}[-2]\br\ar[r]^-{+1}\ar[d]^\bwr&{}\\
M_{\{0\}}[-2]\ar[r]^-2&M_{\{0\}}[-2]\ar[r]&(M/2M)_{\{0\}}[-2]\ar[r]^-{+1}&{}
}
$$
Here $\vphi_x$ is the vanishing cycle functor.
\enEx

\end{document}